\documentclass[11pt,twoside,notitlepage]{article}
\usepackage{a4}
\usepackage[dvips]{graphicx}
\usepackage{amsmath}
\usepackage{amsthm}
\usepackage{amssymb}
\usepackage{amsfonts}
\usepackage{mathrsfs}
\usepackage[all]{xy}
\numberwithin{equation}{section}

\theoremstyle{plain}
\newtheorem{theorem}[subsubsection]{Theorem}

\newtheorem{lemma}[subsubsection]{Lemma}

\newtheorem{prop}[subsubsection]{Proposition}

\newtheorem{question}[subsubsection]{Question}

\theoremstyle{definition}
\newtheorem{definition}[subsubsection]{Definition}

\theoremstyle{remark}
\newtheorem{example}[subsubsection]{Example}
\newtheorem{remark}[subsubsection]{Remark}

\newcommand{\Hc}{\mathcal{H}}

\newcommand{\Oc}{\mathcal{O}}

\newcommand{\Dc}{\mathcal{D}}
\newcommand{\Bc}{\mathcal{B}}

\newcommand{\Lc}{\mathcal{L}}
\newcommand{\Mrm}{\mathrm{M}}

\newcommand{\Pc}{\mathcal{P}}
\newcommand{\Ec}{\mathcal{E}}

\newcommand{\Rc}{\mathcal{R}}

\newcommand{\Uc}{\mathcal{U}}

\newcommand{\Ac}{\mathcal{A}}
\newcommand{\Sf}{\mathfrak{S}}
\newcommand{\Yf}{\mathfrak{Y}}
\newcommand{\Xf}{\mathfrak{X}}
\newcommand{\Zf}{\mathfrak{Z}}

\newcommand{\A}{\mathbb{A}}
\newcommand{\C}{\mathbb{C}}
\newcommand{\R}{\mathbb{R}}
\newcommand{\Q}{{\mathbb{Q}}}

\newcommand{\qb}{{\overline{\Q}}}
\newcommand{\Z}{\mathbb{Z}}

\newcommand{\Zh}{\widehat{\Z}}
\newcommand{\OFh}{{\widehat{\Oc}_F}}
\newcommand{\Nm}{\mathrm{Nm}}
\newcommand{\N}{\mathbb{N}}
\newcommand{\pfk}{\mathfrak{p}}

\newcommand{\mfk}{\mathfrak{m}}

\newcommand{\G}{\mathbb{G}}
\newcommand{\Sb}{\mathbb{S}}
\newcommand{\Hb}{\mathbb{H}}

\newcommand{\Ker}{\mathrm{Ker}}

\newcommand{\Ham}{H}

\newcommand{\Sym}{\mathrm{Sym}}

\newcommand{\Out}{\mathrm{Out}}

\newcommand{\Inn}{\mathrm{Inn}}

\newcommand{\Aut}{\mathrm{Aut}}

\newcommand{\Hom}{\mathrm{Hom}}
\newcommand{\ospace}{\textsc{(OSpace)}}

\newcommand{\Top}{\textsc{(Top)}}

\newcommand{\card}{\mathrm{card}}
\newcommand{\KMS}{\mathrm{KMS}_\beta}

\newcommand{\im}{\mathrm{im}}
\newcommand{\End}{\mathrm{End}}

\newcommand{\princ}{\mathrm{princ}}
\newcommand{\stacks}{\textsc{(Stacks)}}

\newcommand{\id}{\mathrm{id}}

\newcommand{\Id}{\mathrm{Id}}

\newcommand{\KMon}{{K_M}}
\newcommand{\KMonA}{{K^M_\A}}

\newcommand{\Tr}{\mathrm{Tr}}
\newcommand{\Trace}{\mathrm{Trace}}
\newcommand{\Lie}{\mathrm{Lie}}

\newcommand{\limproj}{{\underset{\longleftarrow}{\mathrm{lim}}\,}}

\newcommand{\GL}{\mathrm{GL}}
\newcommand{\SL}{\mathrm{SL}}

\newcommand{\MSp}{\mathrm{MSp}}
\newcommand{\GSp}{\mathrm{GSp}}

\newcommand{\Shf}{\mathfrak{Sh}}
\newcommand{\Sh}{\mathrm{Sh}}

\newcommand{\Res}{\mathrm{Res}}

\newcommand{\innerprod}[2]{\langle#1,#2\rangle}
\newcommand{\norm}[1]{\|#1\|}

\title{Bost-Connes-Marcolli systems\\
for Shimura varieties.\\
I. Definitions and formal analytic properties.}

\author{Eugene Ha and Fr\'ed\'eric Paugam}
\date{July, 2005}

\begin{document}
\maketitle

\abstract{We construct a Quantum Statistical Mechanical system
$(A,\sigma_t)$ analogous to the Bost-Connes-Marcolli
system of \cite{Connes-Marcolli-I} in the case of Shimura varieties.
Along the way, we define a new Bost-Connes system for number fields
which has the ``correct'' symmetries and the ``correct'' partition function.
We give a formalism that applies to general Shimura data $(G,X)$.
The object of this series of papers is to show that these systems have
phase transitions and spontaneous symmetry breaking,
and to classify their KMS states, at least for low temperature.}

\newpage
\tableofcontents
\newpage


\section{Introduction}

A few years ago, Bost and Connes  \cite{Bost-Connes}
discovered a surprising relationship between the class field theory of $\Q$
and quantum statistical mechanics.

Mathematically, a quantum statistical mechanical system consists of a pair
$(\Ac,\sigma_t)$, where $\Ac$ is a C*-algebra and $\sigma_t$ is a
one-parameter group of automorphisms of $\Ac$; physically, $\Ac$ is the algebra of
observables and $\sigma_t$ is the time evolution of the physical system.
The physical states of the system are given by certain linear functionals
on $\Ac$.

The analogy between classical and quantum statistical mechanics can be described
by the following array:\\

\noindent
\begin{tabular}{|c|c|c|}
\hline
    &  \sc Classical        & \sc Quantum\\
\hline
Observables & $a\in C^\infty(X)$    & $a\in \Ac$\\
    & $(X,\omega)$ $2n$-dim symplectic manifold & $\Ac$: C*-algebra, $a$=$a$*\\
    &  (phase space) &\\\hline
Bracket & Poisson bracket       & Commutator\\
        & $\{a_1,a_2\}=\omega(\xi_{a_1},\xi_{a_2})$ & $[a_1,a_2]$\\
    & with $da_i+\omega(\xi_{a_i},\_)=0$ &\\\hline
Hamiltonian & $H\colon X\to \R$     & $H$ unbounded selfadjoint on $\Hc$\\
    &               & Representation $\pi:\Ac\to \Bc(\Hc)$\\\hline
Time    & Solution of           & $\sigma:\R\to \Aut(\Ac)$ \\
evolution& $\{H,a\}(x)=({\frac{d}{dt}})_{t=0}a(\sigma_t(x))$ &
$e^{itH}\pi(a)e^{-itH}=\pi(\sigma_t(a))$\\\hline

States  & Probability measure $\mu$ on $X$  & Linear functional of norm $1$\\
    & $\Phi(a)=\int_X a\mathrm{d}\mu$       & $\Phi:\Ac\to \C$\\\hline
Partition & $\zeta(\beta)=\int_X e^{-\beta H}d\Omega$ & $\zeta(\beta)=\Tr(e^{-\beta H})$\\
function & with $\Omega=\omega^{\wedge n}$ the volume form&\\\hline
Equilibrium & Canonical ensemble &
KMS condition:\\
States & $d\mu=\frac{e^{-\beta H}d\Omega}{\zeta(\beta)}$ & $\Phi(ab)=\Phi(\sigma_{i\beta}(b)a)$\\
&  & example:    $\Phi(a)=\frac{\Tr(ae^{-\beta H})}{\zeta(\beta)}$\\\hline
\end{tabular}\\

The statistical content means that one singles out the equilibrium states
at a given temperature $T=1/\beta$ on $\Ac$,
and these are characterized by the so called $\KMS$ condition.
The set of these equilibrium states may have symmetries.
Changing the temperature of a system can produce a phase transition phenomenon
with spontaneous symmetry breaking, meaning that the symmetry changes radically with an
arbitrary small change of temperature.
For example, the formation at zero temperature of a snowflake from water
is a phase transition, for which we can observe a symmetry
breaking phenomenon: a snowflake has much more symmetry
(it has crystal structure) than a drop of water (which consists
of a random collection of molecules).

The Bost-Connes system $(\Ac,\sigma_t)$
also exhibits a phase transition phenomenon with symmetry
breaking at $\beta=1$. For $\beta<1$, i.e., at high temperature,
there is enough disorder so that the symmetry is trivial. For $\beta>1$,
the set of equilibrium states ``freezes'' and has as symmetry
the Galois group of the maximal abelian extension of $\Q$.
Bost and Connes also defined explicitly a rational
subalgebra of $A_\Q\subset\Ac$ such that the evaluation of KMS states on $A_\Q$,
at small temperature, generate $\Q^{\mathrm{ab}}$.
This system is related to $\GL_{1,\Q}$.

Much more recently, Connes and Marcolli \cite{Connes-Marcolli-I}
defined an analogous system for $\GL_{2,\Q}$, and overcame extreme
technicalities to give in this case
a meaning to all prominent features of the Bost-Connes system
(symmetries, rational subalgebra, zeta function as partition function,
relation to the Galois group of the modular field and its modular reciprocity law).
One of the key points in their new approach is that their system
is related to the study of the ``noncommutative space'' of $\Q$-lattices
up to commensurability:
$$\GL_2(\Q)\backslash \Mrm_2(\A_f)\times \Hb^\pm;$$
and that the set of $\KMS$ states at small temperature is in natural bijection
with the Shimura variety
$$\Sh(\GL_2,\Hb^\pm)=\GL_2(\Q)\backslash \GL_2(\A_f)\times \Hb^\pm.$$
The C*-algebra corresponding to the ``noncommutative space''
of $\Q$-lattices up to commensurability is a groupoid C*-algebra.
This also gives a nice explanation for the origin of the Bost-Connes system.

We choose one direction of generalizing this work of Connes and Marcolli
by replacing their basic Shimura datum $(\GL_2,\Hb^\pm)$ by a general
Shimura datum $(G,X)$.
In order to deal with the technical issue of defining the partition function,
the construction of the Connes-Marcolli system involves a groupoid $\Uc$,
corresponding to the commensurability relation on $\Q$-lattices,
and the quotient of $\Uc$
by the arithmetic subgroup $\SL_2(\Z)\subset \GL_2(\Q)$. We start by defining
an algebra in more adelic terms, meaning that we use a quotient by the
compact open subgroup $\GL_2(\Zh)\subset \GL_2(\A_f)$. The motivation for
this construction comes from the fact that for number fields, one wants
the partition function to be the Dedekind zeta function,
and this is easier to obtain in the adelic language (as pointed out by
Paula Cohen \cite{Cohen}).

The Connes-Marcolli algebra is not exactly a groupoid algebra, because
the quotient of the groupoid $\Uc$ by $\SL_2(\Z)$ is no longer
a groupoid, since $\SL_2(\Z)$ does not act freely on $\Hb$. In fact, if
we use the stacky quotient, then this is a groupoid, but
one cannot define an associated convolution C*-algebra because there is no
good notion of functions on stacks. There are two solutions to this problem,
corresponding to two resolutions of the stack's singularities.
The first is to choose a smaller $\Gamma\subset \SL_2(\Z)$ that acts
freely on $\Hb$. This gives a finite resolution of the stack singularities.
However, this first method works only for classical Shimura varieties,
which does not include the case of a general number field.
The second solution is to identify functions on
$\SL_2(\Z)\backslash \Hb=\GL_2(\Z)\backslash \GL_2(\R)/\C^\times$ to
functions on $\GL_2(\Z)\backslash \GL_2(\R)$ (which is another infinite
resolution of stack singularities) which are invariant for the scaling
action of $\C^\times$. This allows one to define a convolution algebra.
This second method was the one chosen by Connes and Marcolli.

The role of $\Mrm_{2,\Q}$ in the $\GL_{2,\Q}$ case is, in the case of
general Shimura data $(G,X)$, played by a multiplicative semigroup $\Mrm$
such that $\Mrm^\times=G$.
We learned a lot about such semigroups from N. Ramachandran and L. Lafforgue.
Their main properties are given in the appendix.

This article describes the first steps in our work on these Bost-Connes-Marcolli
systems for general Shimura data.

We solve along the way the problem of defining a Bost-Connes system for general
number fields, which has the Dedekind zeta function as partition function
\emph{and} the group of connected components of the idele class group as symmetry
group.
For imaginary quadratic fields, this problem was very recently solved by
Connes-Marcolli-Ramachandran \cite{Connes-Marcolli-Ramachandran}.
Previous works were either restricted to class number one, or did not have the right
symmetry, or did not have the right partition function. All these interesting
works however gradually improved and simplified the techniques involved in the
study of Bost-Connes systems,
and we also use methods from this litterature to prove some of our results.
A generalization of the work \cite{Bost-Connes} to the case of arbitrary
global fields was proposed by Harari and Leichtnam \cite{Harari-Leichtnam}.
A Hecke algebra construction using semi-group crossed products was
proposed in the number field case by Arledge-Laca-Raeburn
\cite{Laca-Raeburn-Arledge}, see also \cite{Laca-Raeburn} and
\cite{Laca}.
Van Frankenhuijsen and Laca \cite{Laca-vanFrankenhuijsen} defined
a system with Galois group as symmetry group for totally imaginary fields of class
number one. P. Cohen \cite{Cohen} constructed adelically
a system with the right partition function in the number field case.
For a nice and more complete survey of known results, see \cite{Connes-Marcolli-I},
Section 1.4.

We also study the explicit example of the Hilbert-Blumenthal modular varieties.

\section*{Acknowledgments}
The authors would like to thank the following persons for useful discussions
in the preparation of this paper:
A. Connes, G. Harder, C. Kaiser, M. Laca, L. Lafforgue, V. Lafforgue, B. Noohi,
D. Panov, N. Ramachandran, B. Toen, D. Zagier.
In particular, N. Ramachandran and L. Lafforgue gave us
references and methods to construct enveloping semigroups, that
are key objects in our formalism. B. Toen gave us a simplicial definition of
stack groupoids.

After writing this paper, we learned from V. Lafforgue another construction
of the Bost-Connes algebra for number fields that will certainly be useful
to study finer aspects of those systems in dimension 1.

We thank the Max-Planck-Institut of Bonn for its hospitality and excellent
working conditions during the preparation of this article. We also
thank the Institut des Hautes \'Etudes Scientifiques for hospitality
during the finalization of this article.

We especially thank Matilde Marcolli for suggesting us to
work on the Connes-Marcolli system in the Hilbert modular case, for
answering all our questions about the article \cite{Connes-Marcolli-I},
and for freely sharing with us her insights.

\section{Another take on the Bost-Connes system}
Before describing the general setting, we would like to present
the illustrating example of the Bost-Connes system, which illuminates
our general constructions.

\subsection{The classical system}
\label{intro-BC-classic}
The \emph{Bost-Connes groupoid} is given by the partially
defined action of $\Q_+^\times$ on $\Zh$. More precisely,
it is given by
$$Z_{BC}=\{(g,\rho)\in \Q_+^\times\times\Zh\mid g\rho\in \Zh\}.$$
The unit space is $\Zh$, and the source and target maps are given
by $s(g,\rho)=\rho$ and $t(g,\rho)=g\rho$. Composition is given by
$(g_2,\rho_2)\circ (g_1,\rho_1)=(g_2g_1,\rho_1)$ if $g_1\rho_1=\rho_2$.

The Bost-Connes Hecke algebra is simply $\Hc:=C_c(Z_{BC})$, equiped
with the convolution product
$$
(f_1*f_2)(g,\rho)=\sum_{h\in \Q_+^\times,h\rho\in\Zh}
f_1(gh^{-1},h\rho)f_2(h,\rho).
$$

The time evolution on this algebra is given by
\begin{equation}
\label{BC-evolution}
\sigma_t(f)(g,\rho)=g^{it}f(g,\rho).
\end{equation}

For each $\rho_0\in\Zh^\times$,
we define a representation $\pi_0:\Hc\to\Bc(\ell^2(\N^\times))$ by
$$(\pi_0(f)(\xi))(n)=\sum_{h\in \N^\times}f(nh^{-1},h\rho_0)\xi(h).$$

To finish, the Hamiltonian of this system is given by
$$H:\ell^2(\N^\times)\to \ell^2(\N^\times),f(n)\mapsto log(n)f(n).$$

By definition, the partition function
$$\zeta_{BC}(s)=\Tr(e^{-s H})=\sum_{n\in\N^\times} n^{-s}=\zeta(s)$$
is exactly Riemann's zeta function.

\subsection{The same system in adelic terms}
\label{intro-BC-adelic}
We will now give a more complicated description of the Bost-Connes
groupoid, that has the advantage of admitting a direct generalization to other
number fields, whose partition function is the Dedekind zeta function.
This is the first step to be carried out in constructing Bost-Connes
systems for number fields. Moreover, the advantage of this
adelic formulation is that it also makes sense for general
Shimura varieties.

We first remark that the quotient of $\Zh$ by the partially defined action of
$\Q_+^\times$ is the same as the quotient of $\Zh\times \{\pm 1\}$ by
the partially defined action of $\Q^\times$. In fact, this equality of quotient
spaces can be described at the level of groupoids.

Let $U^{\princ}\subset \Q^\times\times \Zh\times \{\pm 1\}$ be the groupoid
of elements $(g,\rho,z)$ such that $g\rho\in \Zh$. This groupoid encodes
the partially defined action of $\Q^\times$ on $\Zh\times\{\pm 1\}$.

Now consider the quotient $Z^{\princ}$ of $U^\princ$
by the action of $(\Z^\times)^2=\{\pm 1\}^2$ given by
$$
(\gamma_1,\gamma_2).(g,\rho,z):=
(\gamma_1g\gamma_2^{-1},\gamma_2\rho,\gamma_2z).
$$
This is also a groupoid.

There is a natural morphism of groupoids
$Z_{BC}\to Z^\princ$ given by $(r,\rho)\mapsto (r,\rho,1)$,
which is in fact an isomorphism (cf. \ref{BC+}).

This new descrition of the Bost-Connes groupoid is nicer because
it clearly relates the Bost-Connes system with the pair
$(\G_m,\{\pm 1\})$, which is called the \emph{multiplicative Shimura datum}.

To make this connexion clearer, it is natural to seek a fully adelic description
of the Bost-Connes groupoid. This is because Shimura varieties are
defined adelically. The adelic framework also facilitates the definition of
Bost-Connes systems for number fields with Dedekind zeta function as partition function.

Recall that $\A_f:=\Zh\otimes_\Z\Q$. The strong approximation property
for the multiplicative group $\G_m$ (which in this case is simply
the chinese reminder theorem) tells us that
$$\A_f^\times=\Q_+^\times.\Zh^\times=\Zh^\times.\Q_+^\times.$$
We will now denote
$$\Sh(\G_m,\{\pm 1\}):=\Q^\times\backslash \{\pm 1\}\times \A_f^\times$$
and
$$Y:=\Zh\times \Sh(\G_m,\{\pm 1\}).$$

Consider the partially defined action of $\A_f^\times$ on $Y$ given
by
$$g.(\rho,[z,l]):=(g\rho,[z,lg^{-1}])$$
and let
$$U\subset \A_f^\times\times Y$$
be the corresponding groupoid of elements
$(g,y)$ such that $gy\in Y$.

Now consider the quotient $Z$ of $U$ by the action of $(\Zh^\times)^2$
given by
$$
(\gamma_1,\gamma_2).(g,y):=
(\gamma_1g\gamma_2^{-1},\gamma_2y).
$$

The strong approximation theorem for $\G_m$ implies that the natural map
$$
\begin{array}{ccc}
\Q^\times \times \Zh \times \{\pm 1\} & \to &
\A_f^\times \times \Zh \times \Sh(\G_m,\{\pm 1\})\\
(g,\rho,z) & \to & (g,\rho,[z,1])
\end{array}
$$
induces an isomorphism of groupoids
$Z^\princ\to Z$ (cf. \ref{principal-full}).

The Bost-Connes algebra can thus be described as the algebra
$\Hc=C_c(Z)$ with the convolution product
$$
(f_1*f_2)(g,y)=\sum_{h\in \Zh^\times\backslash\A_f^\times,hy\in Y}
f_1(gh^{-1},hy)f_2(h,y).
$$

Using the isomorphism $d:\Zh^\times\backslash\A_f^\times\to \Q_+^\times$,
we can define the time evolution by
$$\sigma_t(f)(g,y)=d(g)^{it}f(g,y).$$
This is exactly the time evolution we defined in \ref{BC-evolution}.

Let $\Zh^\natural:=\A_f^\times\cap\Zh$ and $\Z^\natural:=\Z-\{0\}$.
The strong approximation theorem gives us that
$$\Zh^\times\backslash \Zh^\natural\cong \Z^\times\backslash \Z^\natural
\cong \N^\times.$$
Let $\Hc_0:=\ell^2(\Zh^\times\backslash \Zh^\natural)\cong \ell^2(\N^\times).$
For each $\rho_0\in\Zh^\times$,
we define a representation $\pi_0:~\Hc\to\Bc(\Hc_0)$ by
$$(\pi_0(f)(\xi))(n)=\sum_{h\in \N^\times}f(nh^{-1},h\rho_0)\xi(h).$$

To finish, the Hamiltonian of this system is given by
$$H:\Hc_0\to \Hc_0,f(n)\mapsto log(d(n))f(n).$$

This adelic system is perfectly identical to the original Bost-Connes
system.
We will however see in the sequel that essentially the same definitions of
algebra, time evolution, representations and Hamiltonian now work
for general Shimura data.

\section{Background material}
In this paper we draw upon the theory of Shimura varieties and operator
algebras.  Since these fields have traditionally had little to do with
each other, we review for the convenience of the reader some of the
basic (well-known) results that we shall need.  This also allows us to
establish notation.  We stress that our definition of a Shimura variety
is a slight variation on the usual one given by Deligne \cite{De4}, 2.1.

\subsection{Shimura varieties}
If $G$ is a reductive group over $\Q$, $G(\R)^+$ will denote the connected
component of identity in the real Lie groups of its real points and
$G(\Q)^+:=G(\Q)\cap G(\R)^+$.
First recall briefly the definition of Shimura data.
We will use a mix of Deligne's definition (see \cite{De4}, 2.1) and Pink's
definition (see \cite{Pink}, 2.1).
Let $\Sb:=\Res_{\C/\R}\G_m$.
\begin{definition}
\label{defshimuradatum}
A \emph{Shimura datum} is a triple $(G,X,h)$, with $G$ a connected
reductive group over $\Q$,
$X$ a left homogeneous space under $G(\R)$ and $h:X\to \Hom(\Sb,G_\R)$ a
$G(\R)$-equivariant map \footnote{for the natural conjugation action of $G(\R)$
on $\Hom(\Sb,G_\R)$} with finite fibres such that:
\begin{enumerate}
\item For $h_x\in h(X)$, $\Lie(G_\R)$ is of type $\{(-1,1),(0,0),(1,-1)\}$;
\item The involution $\mathrm{int}\,h_x(i)$ is a Cartan involution of the
adjoint group $G^{ad}_\R$;
\item The adjoint group has no factor $G'$ defined over $\Q$ on which the
projection of $h_x$ is trivial.
\end{enumerate}
A Shimura datum is called \emph{classical} if it moreover fulfils the
axiom
\begin{enumerate}
\item[4.] Let $Z_0(G)$ be the maximal split subtorus of the center of $G$;
then  $\mathrm{int}\,h_x(i)$ is a Cartan involution of $G/Z_0(G)$.
\end{enumerate}
\end{definition}

\begin{example}
Let $F$ be a number field, $T=\Res_{F/\Q}\G_{m,F}$ and $X_F=T(\R)/T(\R)^+$.
We have $F\otimes_\Q\R\cong \C^i\times \R^j$. We put on $F\otimes_\Q\R$
the Hodge structure that is trivial on $\R^j$ and given by the choice of
a complex structure on $\C^i$ (among the $2^i$ possibilities).
This gives a morphism $h_1:\Sb\to T_\R$. The triple
$(\Res_{F/\Q}\G_{m,F},X_F,h_1)$ is called the
\emph{multiplicative Shimura datum of the field $F$}.
This Shimura datum is classical if and only if $F=\Q$ of $F$ is imaginary quadratic.
\end{example}

We will often denote a Shimura data just by a couple $(G,X)$ when the morphism
$h$ is clear from the situation.

\begin{example}
Let $h:\Sb\to \GL_{2,\R}$ be the morphism given by
$h(a+ib)=\left(\begin{smallmatrix}
a & b\\
-b & a
\end{smallmatrix}\right)$. Let $\Hb^\pm$ be the $\GL_2(\R)$-conjugacy
class of $h$. It identifies with the Poincar\'e double half plane with action
of $\GL_2(\R)$ by homographies. Then $(\GL_2,\Hb^\pm)$ is called the
\emph{modular Shimura datum}.
\end{example}

\begin{definition}
Let $(G,X)$ be a Shimura datum.
Let $K\subset G(\A_f)$ be a compact open subgroup.
The \emph{level $K$ Shimura variety} is
$$\Sh_K(G,X):=G(\Q)\backslash X\times G(\A_f)/K$$
and the \emph{Shimura variety} is the projective limit
$$\Sh(G,X):=\limproj_K \Sh_K(G,X)$$
over all compact open subgroups $K\subset G(\A_f)$.
\end{definition}

A \emph{topological stack} will be for us a stack on the site $\Top$
of topological spaces with usual open coverings, i.e., a category
fibered in groupoids fulfilling some descent condition. See appendix
\ref{stack-groupoids} for more details.

We first remark that, from the modular viewpoint, it is more natural to
study the \emph{level K Shimura stack}, given by the topological stacky quotient
$$\Shf_K(G,X):=[G(\Q)\backslash X\times G(\A_f)/K],$$
and the \emph{Shimura stack}, given by the 2-projective limit
$$\Shf(G,X):=\limproj_K \Shf_K(G,X).$$

In the case of the multiplicative Shimura datum of a number field, i.e.,
$G=\Res_{F/\Q}\G_{m,F}$, the level $K$ Shimura stack can have infinite isotropy
groups given by $G(\Q)^+\cap K$. These isotropy groups are given by generalized
congruence relations on the group of units $\Oc_F^\times$.
We will have to keep track of (some of) these isotropy groups in the
case of non-classical Shimura data.

There is a natural right action of $G(\A_f)$ on $\Sh(G,X)$ given,
for each $g\in G(\A_f)$ and $K\subset G(\A_f)$ compact open, by an isomorphism
$$
\begin{array}{cccc}
(.g): & \Shf_K(G,X) & \to   & \Shf_{g^{-1}Kg}(G,X)\\
      & [x,l]      &\mapsto& [x,lg].
\end{array}
$$

Under the hypothesis that $(G,X)$ is classical
(see also \cite{De4}, 2.1.1.4, 2.1.1.5), there is an easier
description of the Shimura variety (see \cite{De4}, Corollaire 2.1.11):
$$\Sh(G,X)\cong G(\Q)\backslash X\times G(\A_f)$$
and the Shimura stacks $\Shf_K(G,X)$ are in fact algebraic stacks over $\C$.

Unfortunately, these hypothesis are not always fulfilled in the case of the
multiplicative Shimura datum of a general number field $F$. In fact, the quotient
$G(\Q)\backslash X\times G(\A_f)$ is not always Hausdorff in this case.

For example, if $F=\Q(\sqrt{2})$,
$$\Sh(\Res_{F/\Q}\G_{m,F},X_F)\ncong F^\times\backslash X_F\times \A_{f,F}^\times$$
(this is essentially due to the fact that the group of units $\Oc_F^\times$
is infinite).

Points in a Shimura variety will be denoted by pairs $[z,l]$. If the Shimura datum
is classical, this means that $z\in X$ and $l\in G(\A_f)$. Otherwise,
$[z,l]=[z_K,l_K]_{K\subset G(\A_f)}$ is a family of points in $\Sh_K(G,X)$
indexed by the set of compact open subgroups in $G(\A_f)$.

\begin{definition}
Let $(G,X)$ be a Shimura datum.
A compact open subgroup $K\subset G(\A_f)$ is called \emph{neat} if it acts
freely on $G(\Q)\backslash X\times G(\A_f)$.
\end{definition}

We would like to be able to define \emph{natural} algebras of continuous ``functions'' on the finite Shimura varieties in play. In order to do that,
we have to resolve their stack singularities.

We remark that if $K$ is neat, then the quotient analytic stack
$$\Shf_K(G,X)=[G(\Q)\backslash X\times G(\A_f)/K]$$
is a usual analytic space, but otherwise it is worthwhile from the moduli viewpoint
to keep track of the nontrivial stack structure. For classical Shimura data, one can resolve
the stack singularities by choosing a smaller compact open subgroup $K'\subset G(\A_f)$
that acts freely. This is what we will usually do in order to be able to define
continuous ``functions'' on the stack $\Shf_K(G,X)$.

However, this finite resolution of the stack singularities is usually not possible
for nonclassical Shimura data, as
we can see on the following example. Let $F=\Q(\sqrt{2})$ and
$(\Res_{F/\Q}\G_{m,F},X_F)$ be
the corresponding Shimura datum (here $X_F\cong \{\pm 1\}^2$).
Let $K=\OFh^\times$ and consider the stack
$\Shf_K(\Res_{F/\Q}\G_{m,F},X_F)$. Its coarse quotient is the ideal class group of $F$, i.e.,
the trivial group $\{1\}$. Since $F$ has class number one, this coarse quotient
can also be described as $\Oc_F^\times\backslash X_F$.
In this case, $\Oc_F^\times$ is infinite, so that
we can not choose a smaller $K'\subset K$ that acts freely on
$F^\times\backslash X_F\times \A_{f,F}^\times$.
In fact, the unit group is finite if and only if $F=\Q$ or $F$ is imaginary
quadratic, i.e., if and only if
$(\Res_{F/\Q}\G_{m,F},X_F)$ is classical in our language.

If we want to resolve the stack singularities, we can use the quotient map
$$F^\times\backslash \A_F^\times/K\to \Shf_K(\Res_{F/\Q}\G_{m,F},X_F)$$
for the scaling action of the connected component of identity $D_F$ in the
full idele class group $C_F:=F^\times\backslash\A_F^\times$.

\begin{remark}
From the viewpoint of moduli spaces, it is important that the coarse space
$\Sh_K(\Res_{F/\Q}\G_{m,F},X_F)$, i.e.,
the big ideal class group, be replaced by the corresponding
group stack with infinite stabilizers (given by groups of units with congruence
conditions):
$$\Shf_K(\Res_{F/\Q}\G_{m,F},X_F).$$
This ``equivariant viewpoint'' of the finite level Shimura variety could also be important
to understand geometrically the definition of Stark's zeta functions,
and also for the understanding of Manin's real multiplication program \cite{Manin3}.
\end{remark}

\subsection{C*-algebras and quantum statistical mechanics}
We review here basic definitions from the theory of C*-algebras,
emphasising those parts relevant to quantum statistical mechanics. Good
references for the material in this section are \cite{bratteli-robinson-I} and
\cite{bratteli-robinson-II}. For an overview of the grand physical
picture, see \cite{Haag}.

\begin{definition}
\label{C*-condition}
A \emph{C*-algebra} is a (not necessarily unital) complex algebra~$A$
endowed with a conjugate-linear involutive anti-automorphism ${}^*\colon
A\to A$, and a norm~$\norm{\cdot}$, satisfying the following conditions:
For every $a$, $b\in A$ we have
\begin{enumerate}
\item $A$ is complete with respect to the norm, and
$\norm{ab}\le\norm{a}\norm{b}$ (i.e., $A$ is a \emph{Banach algebra});
and
\item $\norm{a^*a}=\norm{a}^2$, the crucial
\emph{C*-condition}.
\end{enumerate}
\end{definition}

Actually a C*-algebra is not as abstract as it may seem, because every
C*-algebra can be realized as a norm-closed sub-*-algebra of the algebra of
bounded operators on a Hilbert space (Theorem of Gelfand-Naimark
\cite{bratteli-robinson-I}, Theorem~2.1.10), and every
such subalgebra is a C*-algebra.
\medskip

The operator algebraic formulation of quantum statistical mechanics (see
the introduction to \cite{bratteli-robinson-I})
consists of a C*-algebra~$A$ together with a 1-parameter group of
automorphism~$\sigma_t\colon A\to A$, which is continuous in the sense
that $t\mapsto\sigma_t(a)$ is continous for every~$a\in A$.  The
algebra~$A$ is then the algebra of quantum observables, while $\sigma_t$
is the time evolution.  The pair~$(A,\sigma_t)$ is an example of a
\emph{C*-dynamical system}.  The \emph{states} of the C*-algebra~$A$ are
the continuous complex-linear functionals~$\Phi$ of norm~$1$ which are
positive, i.e., $\Phi(a^*a)\ge0$ for every~$a\in A$.  The number
$\Phi(a)$ is then the expectation value of the observable $a$ in the
physical state $\Phi$.

To regard the pair~$(A,\sigma_t)$ as a statistical mechanical system we
need an appropriate notion of an ``equilibrium state'' at
temperature~$T=1/\beta$.  This is provided by the KMS condition.

\begin{definition}
The \emph{KMS-$\beta$ condition} ($0<\beta<\infty$) for a state~$\Phi$
is the condition: For every pair of elements $a$, $b\in A$, there is a
complex-valued function~$F$ on the closed
strip $\Omega=\{\,z\in\C\mid0\le\im z\le\beta\,\}$ such that
\[
F(t)=\Phi\bigl(a\sigma_t(b)\bigr),\quad
F(t+i\beta)=\Phi\bigl(\sigma_t(b)a\bigr);
\]
furthermore, the function~$F$ is required to be bounded and continuous
on~$\Omega$, and analytic on its interior.
\end{definition}

This is the definition one often sees in the literature, although in
practice it is easier to use the following equivalent characterization.

\begin{prop}
\label{KMS-prop}
Let $(A,\sigma_t)$ be a C*-dynamical system, and let $\Phi$ be a state
of~$A$.
\begin{enumerate}
\item (\cite{bratteli-robinson-I}, Corollary~2.5.23) There
is a norm-dense *-subalgebra~$A^{\mathrm{an}}$ of~$A$ such that
for every $a\in A^{\rm an}$, the function $t\mapsto \sigma_t(a)$
can be analytically continued to an entire function.
\item (\cite{bratteli-robinson-I},
Definition~5.3.1 and Corollary~5.3.7) The
state~$\Phi$ is a KMS-$\beta$ state if and only if
\[
\Phi\bigl(a\sigma_{i\beta}(b)\bigr)=\Phi(ba)
\]
for all $a$, $b$ in a norm-dense $\sigma_t$-invariant *-subalgebra
of~$A^{\mathrm{an}}$.
\end{enumerate}
\end{prop}

We now proceed to a description of the structure of the set of
KMS-$\beta$ states.  But before doing so, we need to explain the GNS
construction, which is a method of getting representations of a
C*-algebra from its states; it is a basic, widely used result in the
theory of operator algebras.  We also need to define the notion of a
factor state.  We shall use standard notation: given a Hilbert
space~$\mathcal{H}$, we denote the C*-algebra of all bounded operators
on~$\mathcal{H}$ by $B(\mathcal{H})$, and the inner product
on~$\mathcal{H}$ by $\innerprod{\cdot}{\cdot}$.

\begin{prop}%
[GNS construction; \cite{bratteli-robinson-I}, 2.3.16]
Let $\Phi$ be a state of a C*-algebra~$A$.  Then there is a
triple~$(\mathcal{H}_\Phi,\pi_\Phi,\xi_\Phi)$ consisting of a
representation~$\pi_\Phi$ of~$A$ on a Hilbert space~$\mathcal{H}_\Phi$
and a unit vector~$\xi_\Phi\in\mathcal{H}_\Phi$ such that:
\begin{enumerate}
\item $\Phi(a)=\innerprod{\pi_\Phi(a)\xi_\Phi}{\xi_\Phi}$ for all~$a\in
A$; and
\item The orbit $\pi_\Phi(A)\xi_\Phi$ is norm-dense
in~$B(\mathcal{H}_\Phi)$.
\end{enumerate}
The triple~$(\mathcal{H}_\Phi,\pi_\Phi,\xi_\Phi)$ is unique up to
unitary equivalence.
\end{prop}

The states of particular relevance to the KMS theory are the
\emph{factor states}.  These are the states~$\Phi$ for which the
corresponding GNS representation~$\pi_\Phi$ generates a \emph{factor},
which is to say that the weak closure of $\pi_\Phi(A)$ in
$B(\mathcal{H}_\Phi)$ has centre consisting of the scalar
operators.  (This weak closure is an example of a \emph{Von Neumann
algebra}.)
\medskip

We can now state the main structure theorem for the set of KMS-$\beta$
states.

\begin{prop}
[Structure of KMS states; \cite{bratteli-robinson-II}, Theorem~5.3.30]
\label{KMS-structure}
The set $\mathcal{E}_\beta$ of KMS-$\beta$ states is a convex,
weak*-compact simplex.  The extreme points of~$\mathcal{E}_\beta$ are
precisely those KMS-$\beta$ states that are factor states.
\end{prop}

\section{Abstract Bost-Connes-Marcolli systems}
\label{BCM-general}
The aim of this section is to
define Bost-Connes-Marcolli systems for general Shimura data
$(G,X)$ and study their basic formal properties.
A better understanding of the general setup might
be gained by looking at section \ref{BCM-number-fields} where we specialise to the
case of multiplicative Shimura datum (the case relevant for number fields).

\subsection{BCM data}
\label{BCMdata}
In order to define a generalization of the Connes-Marcolli algebra to general
Shimura data, we want to make clear the separation between algebraic and
level structure data, which is already implicit in the construction of
Connes and Marcolli.

\medskip
\noindent\textit{{Algebraic data}}.\enspace
We first need to consider a semigroup $M$ which plays the role for
a general reductive group $G$ that $\Mrm_{2,\Q}$ plays for $\GL_{2,\Q}$.
\begin{definition}
Let $G$ be reductive group over a field. An \emph{enveloping semigroup}
for $G$ is a multiplicative
semigroup $M$ which is irreducible and normal, and such that $M^\times=G$.
\end{definition}

\begin{definition}
A \emph{BCM datum} is a tuple $\Dc=(G,X,V,M)$ with $(G,X)$ a Shimura datum,
$V$ a faithful representation of $G$ and $M$ an enveloping semigroup
for $G$ contained in $\End(V)$.
\end{definition}

The faithful representation will often be denoted $\phi:G\to \GL(V)$.

\medskip
\noindent\textit{{Level structure data}}.\enspace
Every Shimura datum $(G,X)$ comes implicitly with a family of level structures given by
the family of compact open subgroups $K\subset G(\A_f)$. Connes and Marcolli fixed the
full level structure $\GL_2(\Zh)\subset \GL_2(\A_f)$ as starting datum for their
construction. To avoid the problem they had with stack singularities of their
groupoid, we will fix a neat level structure as part of the datum.

The level structure also plays a role in defining the partition function of our system.
Consideration of maximal level structures then yields standard zeta functions
as partition functions, for example, the Dedekind zeta function of a number field.
A technical requirement in the definition of the partition function is the choice of
a lattice in the representation of $G$, which enables us to define a rational
determinant for the adelic matrices in play.

\begin{definition}
Let $\Dc=(G,X,V,M)$ be a BCM datum. A \emph{level structure on $\Dc$} is a triple
$\Lc=(L,K,\KMon)$, with $L\subset V$ a lattice, $K\subset G(\A_f)$ a compact open
subgroup, and $\KMon\subset M(\A_f)$ a compact open subsemigroup, such that
\begin{itemize}
\item $\KMon$ stabilizes $L\otimes_\Z\Zh$,
\item $\phi(K)$ is contained in $\KMon$.
\end{itemize}
The pair
$(\Dc,\Lc)$ will be called \emph{a BCM pair}.
\end{definition}

We can summarize the relation between $L$, $K$ and $\KMon$ by the following diagram:
$$
\xymatrix@+10pt{
K\vphantom{\int}\ar[r]^\phi\ar@{^{(}->}[d] &
\KMon\vphantom{\int}\phantom{\_} \ar@{^{(}->}[r]\ar@{^{(}->}[d] &
\End(L)(\Zh)\vphantom{\int}\phantom{\_}\ar@{^{(}->}[d]\\
G(\A_f)\ar[r]^\phi &
M(\A_f)\phantom{\_}\ar@{^{(}->}[r] &
\End(V)(\A_f)
}
$$

\begin{definition}
The \emph{maximal level structure} $\Lc_0=(L,K_0,K_{M,0})$ associated with
a datum $\Dc=(G,X,V,M)$ and a lattice
$L\subset V$ is defined by setting
$$
\begin{array}{ccc}
K_{M,0} & := & M(\A_f)\cap\End(L\otimes_\Z\Zh),\\
K_0     & := & \phi^{-1}(K_{M,0}^\times).
\end{array}
$$
\end{definition}

\begin{definition}
The level structure $\Lc$ on $\Dc$ is called \emph{fine} if $K$ acts
freely on $G(\Q)\backslash X\times G(\A_f)$.
\end{definition}

The maximal level structure is usually not neat enough to avoid stack
singularity problems in the generalization of the Connes-Marcolli algebra.
This is why we introduce the
additional data of a compact open subgroup $K\subset \KMon$.
For example, for the Connes-Marcolli case, one takes $K=\GL_2(\Zh)$, $\KMon=\Mrm_2(\Zh)$,
but the fact that this choice of $K$ is not neat implies that the groupoid we introduce
in the next section has stack singularities.
Thus we instead choose a smaller $K=K(N)\subset \GL_2(\Zh)$ given
by the kernel of the mod $N$ reduction of matrices.

\medskip
\noindent\textit{Symmetries and zeta function}.\enspace
The symmetries of the Connes-Marcolli system play an important role in its relations
with arithmetic.
The analogous symmetry in our generalization is the following
(which will be justified in Subsection \ref{symmetries-general}).
\begin{definition}
\label{sympm-general}
The semigroup
$\Sym_f(\Dc,\Lc):=\phi^{-1}(\KMon)$ is called the \emph{finite symmetry
semigroup of the BCM pair $(\Dc,\Lc)$}.
We will denote by $\Sym_f^\times(\Dc,\Lc)$ the group of invertible elements
in $\Sym_f(\Dc,\Lc)$.
\end{definition}

We included in $\Lc$ the datum of a lattice in the representation $\phi$
in order to define a determinant map.
\begin{lemma}
\label{lemme-determinant}
The determinant $\det:\GL(L)\to\G_m$ induces a natural map,
$$(\det\circ\,\phi):K\backslash G(\A_f)/K\to \Q^\times_+.$$
The image of $\Sym_f(\Dc,\Lc)$ under this map is contained in $\N^\times$.
\end{lemma}
\begin{proof}
Since $\phi(K)\subset K_M^\times\subset \GL(L)(\Zh)$,
the representation $\phi:G\to \GL(L\otimes_\Z\Q)$ induces a map
$$
\phi:K\backslash G(\A_f)/K\to\GL(L)(\Zh)\backslash\GL(L)(\A_f)/\GL(L)(\Zh).
$$
The determinant map $\GL(L)\to \G_m$ induces a natural map
$$
\det:\GL(L)(\Zh)\backslash\GL(L)(\A_f)/\GL(L)(\Zh)\to
\Zh^\times\backslash\A_f^\times/\Zh^\times\cong \Zh^\times\backslash\A_f^\times\cong
\Z^\times\backslash\Q^\times\cong \Q^\times_+.
$$
The composition $\det\circ\,\phi$ gives us the desired map.
The image of $\Sym_f$ under this map is contained in the image
of $\GL(L)(\A_f)\cap \End(L)(\Zh)$  under the determinant map,
which is exactly $\Zh^\natural:=\A_f^\times\cap \Zh$. The quotient
$\Zh^\times\backslash \Zh^\natural$ is identified with
$\Z^\times\backslash \Z\cong \N^\times\subset \Q^\times$.
\end{proof}

\begin{definition}
\label{zeta-general}
The \emph{zeta function} of the BCM pair $(\Dc,\Lc)$ is the complex valued series
$$
\zeta_{\Dc,\Lc}(\beta):=
\sum_{g\in \Sym_f^\times\backslash\Sym_f}\det(\phi(g))^{-\beta}.$$
The BCM pair $(\Dc,\Lc)$ is called \emph{summable} if
there exists $\beta_0\in \R$ such that
$\zeta_{\Dc,\Lc}(\beta)$ converges in the right plane
$\{\beta\in\C\mid \textrm{Re}(\beta)>\beta_0\}$
and extends to a meromorphic function on the full complex plane.
\end{definition}

\subsection{The BCM groupoid}
\label{BCMgroupoid}
Let $(\Dc,\Lc)=((G,X,V,M),(L,K,\KMon))$ be a BCM pair.
There are left and right actions of $G(\A_f)$ on $M(\A_f)$.

\subsubsection{Definition}
Connes and Marcolli remarked in \cite{Connes-Marcolli-I} that, if we want
to take a quotient of a groupoid by a group action, it is
essential that the action is free on the unit space of the groupoid.
If we take the usual quotient set of a groupoid by an action that is
not free on the unit space, this will not give a groupoid.
We are thus obliged to use unit spaces that are in fact stacks.
Some of them have nice singularities (i.e., those with finite stabilizers).
Others don't, but the language of stacks allows one to work in full
generality without bothering about the freeness of actions in play.

We will denote the stacks by german letters;
the corresponding coarse spaces will be denoted by right letters.

Let
$$Y_{\Dc,\Lc}=\KMon\times \Sh(G,X).$$
We denote points of $Y_{\Dc,\Lc}$ by triples $y=(\rho,[z,l])$ with
$\rho\in \KMon$, $[z,l]\in \Sh(G,X)$.

We want to study the equivalence relation on $Y_{\Dc,\Lc}$
given by the following partially defined action of $G(\A_f)$:
$$g.y=(g\rho,[z,lg^{-1}]),\quad\textrm{where }y=(\rho,[z,l]).$$
This equivalence relation will be called the \emph{commensurability
relation}. This terminology is derived from the notion of commensurability
for $\Q$-lattices, cf. \cite{Connes-Marcolli-I}.

Consider the subspace
$$U_{\Dc,\Lc}\subset G(\A_f)\times Y_{\Dc,\Lc}$$
of pairs $(g,y)$ such that $gy\in Y_{\Dc,\Lc}$, i.e. $g\rho\in \KMon$.

The space $U_{\Dc,\Lc}$ is a groupoid with unit space $Y_{\Dc,\Lc}$.
The source and target maps
$s:U_{\Dc,\Lc}\to Y_{\Dc,\Lc}$ and $t:U_{\Dc,\Lc}\to Y_{\Dc,\Lc}$ are given by
$s(g,y)=y$ and $t(g,y)=gy$.
The composition is given, for $y_1=g_2y_2$,
by $(g_1,y_1)\circ (g_2,y_2)=(g_1g_2,y_2)$.
Notice that the groupoid obtained by
restricting this groupoid
to the $(g,(\rho,[z,l]))$ such that $\rho$ is invertible
is free, i.e., the equality $t(g,y)=s(g,y)$ implies $g=1$.

There is a natural action of $K^2$ on the groupoid $U_{\Dc,\Lc}$, given by
$$(g,y)\mapsto (\gamma_1g\gamma_2^{-1},\gamma_2y),$$
and the induced action on $Y_{\Dc,\Lc}$ is given by
$$y\mapsto \gamma_2y.$$

There are two motivations for quotienting $U_{\Dc,\Lc}$ by this action. The first
one is physical: it is necessary to obtain a reasonable partition function for
our system. The second is moduli theoretic: $U_{\Dc,\Lc}$ is only a pro-analytic
groupoid and the quotient by $K^2$ is fibered over the Shimura variety $\Shf_K(G,X)$
which is an algebraic moduli stack of finite type whose definition
could be made over $\qb$, at least when $(G,X)$ is classical and the Shimura variety
has a canonical model.

Let $\Zf_{\Dc,\Lc}$ be the quotient stack
$[K^2\backslash U_{\Dc,\Lc}]$ and $\Sf_{\Dc,\Lc}$ be the quotient stack
$[K\backslash Y_{\Dc,\Lc}]$.
The natural maps
$$s,t:\Zf_{\Dc,\Lc}\to \Sf_{\Dc,\Lc}$$
define a stack-groupoid structure (see appendix \ref{stack-groupoids}) on
$\Zf_{\Dc,\Lc}$ with unit stack $\Sf_{\Dc,\Lc}$.

\begin{definition}
The stack-groupoid $\Zf_{\Dc,\Lc}$ is called the \emph{Bost-Connes-Marcolli}\footnote{We
will often call it the BCM groupoid, for short.} groupoid.
\end{definition}

Let $Z_{\Dc,\Lc}:=K^2\backslash U_{\Dc,\Lc}$
be the (classical, i.e., coarse) quotient of $U_{\Dc,\Lc}$
by the action of $K^2$.
If $K$ is small enough, i.e., if $K$ acts freely on
$G(\Q)\backslash X\times G(\A_f)$, then $\Zf_{\Dc,\Lc}$
is equal to the classical quotient $Z_{\Dc,\Lc}$, which is a groupoid
in the usual sense,
with units $S=K\backslash Y_{\Dc,\Lc}$. Otherwise, suppose that there
exists  a compact open subgroup $K'\subset K$ that acts freely on
$G(\Q)\backslash X\times G(\A_f)$ and choose on $\Dc$ the level structure
$\Lc'=(L,K',\KMon)$. The stack $\Zf_{\Dc,\Lc'}$ is a usual topological space
that is a finite covering of the coarse space $Z_{\Dc,\Lc}$ and such that
the stack $\Zf_{\Dc,\Lc}$ is the stacky quotient of $Z_{\Dc,\Lc'}$ by the projection
equivalence relation to $Z_{\Dc,\Lc}$.

The reader who prefers to work
with usual analytic spaces will thus suppose that $K$ is small enough,
but as we remarked before, our basic examples (number fields) do not fulfil this hypothesis.
We have also to recall that for nonclassical Shimura data $(G,X)$ in the sense of
definition \ref{defshimuradatum}, there exists no such small enough
$K\subset G(\A_f)$. This is essentially due to the fact that the ``unit group''
$C(\Q)\cap K$ (where $C$ denotes the center of $G$) can be infinite.

\subsubsection{The commensurability class map}
Recall that $\phi$ is the representation of $G$, which we view as the natural
inclusion $G\to M$.

\medskip
\noindent\textit{{For classical Shimura data}}.\enspace
We want to give an explicit description of the quotient of $Y_{\Dc,\Lc}$
by the commensurability equivalence relation, in the case where $(G,X)$ is
classical, i.e., when
$$\Sh(G,X)=G(\Q)\backslash X\times G(\A_f).$$

Let $\KMonA=\phi(G)(\A_f).\KMon\subset M(\A_f)$.
There is a natural surjective map of sets
$$\pi:Y_{\Dc,\Lc}\to G(\Q)\backslash X\times \KMonA$$
given by
$\pi(\rho,[z,l])=[z,l\rho]$.

Let $Y_{\Dc,\Lc}^\times=K_M^\times\times (G(\Q)\backslash X\times G(\A_f))$
be the invertible part of $Y_{\Dc,\Lc}$ and let
$Z_{\Dc,\Lc}^\times\subset \Zf_{\Dc,\Lc}$ be
the corresponding subspace (which is a groupoid in the usual sense because $K$ acts freely on $K_M^\times$); that is, $Z^\times_{\Dc,\Lc}$ is defined just as
$\Zf_{\Dc,\Lc}$ is, but with $Y^\times_{\Dc,\Lc}$ in place of $Y_{\Dc,\Lc}$.
Let $S_{\Dc,\Lc}^\times:=K\backslash Y_{\Dc,\Lc}^\times$ be the unit space of $Z_{\Dc,\Lc}^\times$.
Since $K_M^\times\subset M^\times(\A_f)=\phi(G(\A_f))$, the map $\pi$
induces a natural map
$$
\begin{array}{cccc}
\pi^\times: & Y_{\Dc,\Lc}^\times     & \to   & G(\Q)\backslash X\times G(\A_f)\\
     & (\rho,[z,l]) &\mapsto& [z,l\phi^{-1}(\rho)],
\end{array}
$$
which is complex analytic (for the natural analytic structures induced by
the complex structure on $X$) and surjective.
Both $\pi$ and $\pi^\times$ factor through the quotient of their sources
by the left action of $K$. We will continue to denote this factorisation by
$\pi$ and $\pi^\times$.

\begin{definition}
The maps $\pi$ and $\pi^\times$ are called the \emph{commensurability class maps}.
\end{definition}

The last definition is justified by the following lemma. The notion
of coarse quotient can be found in Definition \ref{coarse-quotient}.
\begin{lemma}
\label{quotient-of-groupoid}
The maps $\pi$ and $\pi^\times$ are in fact
the coarse quotient maps for the groupoids $\Zf_{\Dc,\Lc}$ and
$Z_{\Dc,\Lc}^\times$ acting on their
unit spaces $\Sf_{\Dc,\Lc}$ and $S_{\Dc,\Lc}^\times$.
\end{lemma}
\begin{proof}
If $(g,\rho,[z,l])\in U_{\Dc,\Lc}$, then
$\pi(g\rho,[z,lg^{-1}])=[z,l\rho]=\pi(\rho,[z,l])$ which proves that $\pi$
factors through
$$|\Sf_{\Dc,\Lc}/\Zf_{\Dc,\Lc}|\to G(\Q)\backslash X\times \KMonA.$$
This surjective map is in fact an isomorphism.
Indeed, if $(\rho,[z,l]),(\rho',[z',l'])\in Y_{\Dc,\Lc}$ have same image under $\pi$,
then there exists $g\in G(\Q)$ such that $gl\rho=l'\rho'$ and $gz=z'$.
We then know that in the quotient space
$|\Sf_{\Dc,\Lc}/\Zf_{\Dc,\Lc}|$,
$$
\begin{array}{ccl}
(\rho,[z,l])    &=&     (l^{-1}g^{-1}gl\rho,[z,l])\\
        &=&     (l^{-1}g^{-1}l'\rho',[z,l])\\
        &\sim&  ({l'}^{-1}gll^{-1}g^{-1}l'\rho',[z,ll^{-1}g^{-1}l'])\\
        &=& (\rho',[z,g^{-1}l'])\\
        &=& (\rho',[gz,l'])\\
        &=& (\rho',[z',l']).
\end{array}
$$
This proves injectivity of $\pi$ and surjectivity was already known.
The argument for $\pi^\times$ is similar.
\end{proof}

\medskip
\noindent\textit{{For commutative Shimura data}}.\enspace
Commutative Shimura data form another family of examples for which we can construct
the commensurability class map in simple terms. The multiplicative datum of a
number field is in this familly.
Thus we now suppose that $\Dc=(G,X,V,M)$ is a BCM datum such that $G$ and
$M$ are commutative, and we let $\Lc$ be a level structure on $\Dc$.
For each $K',K\subset G(\A_f)$ compact open, there is a natural map
$$\Yf_{K'}:=\KMon\times \Shf_{K'}(G,X)\to [G(\Q)\backslash X\times M(\A_f)/K']$$
given by $(\rho,[z,l])\mapsto [z,l\rho]$.
This map is $K$-equivariant for the trivial action of $K$
on the range because the image of
$k.(\rho,[z,l])=(k\rho,[z,lk^{-1}])$ is equal to the image of $(\rho,[z,l])$.
Recall that $\Sf_{\Dc,\Lc}:=[K\backslash Y_{\Dc,\Lc}]$ and
$S_{\Dc,\Lc}^\times=K\backslash Y_{\Dc,\Lc}^\times$.
If we pass to the limit on $K'\subset G(\A_f)$, and then to the quotient by $K$,
we obtain natural maps
$$
\pi:\Sf_{\Dc,\Lc}\to \limproj_{K'}[G(\Q)\backslash X\times M(\A_f)/K']
$$
and
$$
\pi^\times:S_{\Dc,\Lc}^\times\to \Sh(G,X)
$$
that will be called as before the \emph{commensurability class maps}.

The image of the map $\pi$ is as before the coarse quotient for the action of
the groupoid $\Zf$ on its unit space $\Sf$.

Denote $\Sf_{K'}:=K\backslash (\KMon\times \Shf_{K'}(G,X))$.
We should remark here that in this commutative case, the space
$\Sf_{K'}$ is the unit space of a well defined groupoid
$\Zf_{K'}$ because the $G(\A_f)$ action
on $\Shf_{K'}(G,X)$ is well defined. This shows that
\begin{equation}
\label{equa-limproj-groupoide}
\Zf=\limproj_{K'} \Zf_{K'},
\end{equation}
which will be useful for the description of the symmetries of
Bost-Connes systems for number fields.

\subsection{Defining BCM algebras}
\label{BCMalgebras}
\subsubsection{Functions on BCM stacks?}
Let $\Dc=(G,X,V,M)$ be a BCM datum, $V$ be a representation of $G$, and
let $L_0$ be the associated maximal level structure (Def.~4.1.4).
We would like to define the BCM algebra of $(\Dc,\Lc_0)$ as a groupoid
algebra. Unfortunately, the corresponding groupoid is usually only a stack
and there is no canonical notion of continuous functions on such a space.
More precisely, if a Stack has some nontrivial isotropy group, Connes' philosophy
of noncommutative geometry tells us that the ``algebra of functions'' on it should
include this isotropy information in a nontrivial way, and this algebra depends
on a presentation of the stack.

If $(G,X)$ is classical, there is a very natural way to resolve the
stack singularities of $\Zf_{\Dc,\Lc_0}$ by choosing a neat
level structure $\Lc$, for which the projection map
$$Z_{\Dc,\Lc}\to \Zf_{\Dc,\Lc_0}$$
is such a resolution.
The corresponding convolution algebra of the groupoid $Z_{\Dc,\Lc}$
is a completely natural replacement for the
groupoid algebra of the stack-groupoid $\Zf_{\Dc,\Lc_0}$.

If $(G,X)$ is nonclassical, there is no nice resolution of the stack singularities
of $\Zf_{\Dc,\Lc_0}$. We will thus work with the algebra of functions on
the coarse quotient $Z_{\Dc,\Lc_0}$.
However $Z_{\Dc,\Lc_0}$ is \emph{not} a groupoid,
and so to define a convolution algebra from the function algebra $C_c(Z_{\Dc,\Lc_0})$,
we use the trick used by Connes and Marcolli in \cite{Connes-Marcolli-I}, 1.83,
which consists in introducing $G(\R)$.
Namely, we introduce the groupoid
$$
\Rc_{\Dc,\Lc_0}\subset
K\backslash G(\A_f)\underset{K}{\times} (\KMon\times \limproj_{K'} G(\Q)\backslash G(\A)/K'),
$$
where $K'$ runs over compact open subgroups of $G(\A_f)$;
and identify $C_c(Z_{\Dc,\Lc_0})$ with the subalgebra of $C_c(\Rc_{\Dc,\Lc})$
obtained by composing by the projection map $\Rc_{\Dc,\Lc_0}\to Z_{\Dc,\Lc_0}$.
Since $\Rc_{\Dc,\Lc_0}$ is a groupoid, convolution can be defined on $C_c(Z_{\Dc,\Lc})$.

Remark that this solution, even if not completely satisfactory from the
geometrical viewpoint (because we work on coarse quotients), suffices
(and seems to be necessary) for the physical interpretation, i.e., analysis of KMS states.

\subsubsection{BCM algebras}
Now we give the precise definition of the algebra alluded to in the previous paragraph.
Let $(\Dc,\Lc)=((G,X,V,M)(L,K,\KMon))$ be a BCM pair.

Let
$$\Hc(\Dc,\Lc):=C_c(Z_{\Dc,\Lc})$$
be the algebra of compactly supported continuous functions on $Z_{\Dc,\Lc}$.
As in \cite{Connes-Marcolli-I}, p44, in order to define the convolution of
two functions, we consider functions on $Z_{\Dc,\Lc}$ as functions on
$U_{\Dc,\Lc}$ satisfying the following properties:
$$
f(\gamma g,y)=f(g,y)\quad f(g\gamma,y)=f(g,\gamma y),\quad
\forall \gamma\in K,\quad g\in G(\A_f),\quad y\in Y_{\Dc,\Lc}.
$$

The convolution product on $\Hc(\Dc,\Lc)$ is then defined by the expression
$$
(f_1*f_2)(g,y):=
\sum_{h\in K\backslash G(\A_f),\,hy\in Y_{\Dc,\Lc}}
f_1(gh^{-1},hy)f_2(h,y),
$$
and the adjoint by
$$f^*(g,y):=\overline{f(g^{-1},gy)}.$$

The fact that we consider functions with compact support implies that the
sum defining the convolution product is finite.

\begin{definition}
\label{BCM-algebras}
The algebra $\Hc(\Dc,\Lc)$ (under the convolution product)
is called the \emph{Bost-Connes-Marcolli algebra} of the pair $(\Dc,\Lc)$.
\end{definition}

\begin{remark}
We proved in Lemma \ref{quotient-of-groupoid} that, if $(G,X)$ is classical,
the quotient
of $Y_{\Dc,\Lc}$ by the commensurability equivalence relation
(encoded by the action of the groupoid $Z_{\Dc,\Lc}$) does not depend on
the choice of $K$. This implies that in the classical case, the Morita equivalence class of
$\Hc(\Dc,\Lc)$ is independent of the choice of neat level structure $K$.
More precisely, all these algebras are in fact Morita equivalent to
the algebra corresponding to the ``noncommutative quotient''
$$G(\Q)\backslash X\times K^M_\A,\quad\text{where $K^M_\A=G(\A_f)K_M$.}$$
\end{remark}

\subsection{Time evolution, Hamiltonian and partition function}
\label{tev-hamiltonien-general}
Let $(\Dc,\Lc)=((G,X,V,M),(L,K,\KMon))$ be a BCM pair with neat level.

\begin{definition}
\label{tev-general}
The \emph{time evolution} on $\Hc(\Dc,\Lc)$ is defined by
\begin{equation}
\label{eq:time-evolution}
\sigma_t(f)(g,y)=\det(\phi(g))^{it}f(g,y).
\end{equation}
\end{definition}

Let $y=(\rho,[z,l])$ be in $Y_{\Dc,\Lc}$ and let
$G_y=\{g\in G(\A_f)\mid g\rho\in \KMon\}$. Let
$\Hc_y$ be the Hilbert space $\ell^2(K\backslash G_y)$.

\begin{definition}
The \emph{representation
$\pi_y:\Hc(\Dc,\Lc)\to \Bc(\Hc_y)$ of the Hecke algebra on $\Hc_y$}
is defined by
$$
(\pi_y(f)\xi)(g):=\sum_{h\in K\backslash G_y}
f(gh^{-1},hy)\xi(h),\;\forall g\in G_y,
$$
for $f\in \Hc(\Dc,\Lc)$ and $\xi\in \Hc_y$.
\end{definition}

\begin{lemma}
The representation $\pi_y$ is well defined, i.e., $\pi_y(f)$ is bounded for each
$f\in \Hc(\Dc,\Lc)$.
\end{lemma}
\begin{proof}
For $f\in \Hc(\Dc,\Lc)$,
We want to prove that the norm
$$\|\pi_y(f)\|:=\sup_{\|\xi\|=1}\|\pi_y(f)\xi\|_2$$
is bounded. This follows from the fact that
the functions we consider are with compact support.
More precisely,
denote $Z:=Z_{\Dc,\Lc}$.
Given $f\in \Hc(\Dc,\Lc)=C_c(Z)$, we need to show that there is a bound $C>0$
such that for every pair of vectors $\xi,\eta\in\mathcal{H}_y$ we have
\[
|\innerprod{\pi_y(f)\xi}{\eta}|\le C\norm{\xi}\norm{\eta}.
\]
To this end, we introduce the following notation.  We set
\[
S_y=\{\,[gh^{-1},hy]\in Z\mid g,h\in K\backslash G_y\,\},
\]
and for each $\gamma\in S_y$ we set
\[
R_y(\gamma)=\{\,\gamma'\in Z_y\mid s(\gamma')=t(\gamma)\,\}.
\]
These are discrete sets.  Here we use the usual notation for groupoids,
namely $Z_y=t^{-1}\{y\}$, which we shall identify with
$K\backslash G_y$.

Using the Cauchy-Schwarz inequality, we now get a bound on
$|\innerprod{\pi_y(f)\xi}{\eta}|$ as follows:
\begin{align*}
|\innerprod{\pi_y(f)\xi}{\eta}|
&\le\sum_{\gamma_1\in Z_y}
    \bigl|\bigl(\pi_y(f)\xi\bigr)(\gamma_1)
          \overline{\eta(\gamma_1)}\bigr| \\
&\le\sum_{\gamma_1,\gamma_2\in Z_y}
    \bigl|f(\gamma_1\gamma_2^{-1})\xi(\gamma_2)
                                  \overline{\eta(\gamma_1)}\bigr| \\
&=\sum_{\gamma\in S_y}|f(\gamma)|
  \sum_{\gamma'\in R_y(\gamma)}|\xi(\gamma')\eta(\gamma\gamma')| \\
&\le\sum_{\gamma\in S_y}|f(\gamma)|
   \biggl(\sum_{\gamma'\in R_y(\gamma)}
         |\xi(\gamma')|^2\biggr)^{\frac{1}{2}}
   \biggl(\sum_{\gamma'\in R_y(\gamma)}
         |\eta(\gamma\gamma')|^2\biggr)^{\frac{1}{2}} \\
&\le\norm{\xi}\norm{\eta}\sum_{\gamma\in S_y}|f(\gamma)|.
\end{align*}
Because $f$ has compact support, the sum $\sum_{\gamma\in
S_y}|f(\gamma)|$ is finite, and we thereby get the desired bound.
\end{proof}

Let $K_0=\phi^{-1}(K_M^\times)$.
We view the Hamiltonian as a virtual operator on $\ell^2(K_0\backslash G_y)$.
By this we mean that the Hamiltonian does not depend on the choice of $K$
and there is a minimal space on which it is defined: the space
$\ell^2(K_0\backslash G_y)$.
Consequently, its trace must be computed as a virtual (i.e., equivariant) trace,
i.e., must be divided by $\card(K\backslash K_0)$.
These considerations are related to the
fact that, if $(G,X)$ is classical, we prefer to define BCM algebras using neat
level structures to resolve the stack singularities of $\Zf_{\Dc,\Lc}$.

\begin{prop}
\label{ham-general}
The operator on $\Hc_y$ given by
$$(\Ham_y\xi)(g)=\log\det(\phi(g))\cdot\xi(g)$$
is the \emph{Hamiltonian}, i.e., the infinitesimal generator of the time evolution,
meaning that we have the equality
\begin{equation}
\label{eq:hamiltonian}
\pi_y(\sigma_t(f))=e^{it\Ham_y}\pi_y(f)e^{-it\Ham_y}
\end{equation}
for all $f\in \Hc(\Dc,\Lc)$.
\end{prop}
\begin{proof}
This is just a matter of unwinding the definitions.  Let
$\xi\in\Hc_y$, and let $g\in G_y$.  On the one hand we have
\begin{align*}
\bigl(\pi_y(\sigma_tf)\xi\bigr)(g)
&=\sum_{h\in K\backslash G_y}
  (\sigma_tf)(gh^{-1},hy)\xi(h)\\
&=\sum_{h\in K\backslash G_y}
  \det(\phi(g))^{it}\det(\phi(h))^{-it}f(gh^{-1},hy)\xi(h),
\end{align*}
while on the other hand we have
\begin{align*}
\bigl(e^{it \Ham_y}(\pi_yf)e^{-it\Ham_y}\xi\bigr)(g)
&=\det(\phi(g))^{it}\bigl((\pi_yf)e^{-it\Ham_y}\xi\bigr)(g)\\
&=\det(\phi(g))^{it}\sum_{h\in K\backslash G_y}
  f(gh^{-1},hy)(e^{-it\Ham_y}\xi)(h)\\
&=\det(\phi(g))^{it}\sum_{h\in K\backslash G_y}
  f(gh^{-1},hy)\det(\phi(h))^{-it}\xi(h).
\end{align*}
We thereby obtain the desired equality.
\end{proof}

\begin{definition}
Let $y\in Y_{\Dc,\Lc}$ and $\beta>0$. The \emph{partition function} of the system
$(\Hc(\Dc,\Lc),\sigma_t,\Hc_y,\Ham_y)$,
is
$$\zeta_y(\beta):=\frac{1}{\card(K\backslash K_0)}\Trace(e^{-\beta\Ham_y}).$$
\end{definition}

Let $Y_{\Dc,\Lc}^\times\subset Y_{\Dc,\Lc}$ be the set of invertible $y=(\rho,[z,l])$,
i.e., $\rho\in K_M^\times$.

\begin{prop}
Suppose that $y\in Y_{\Dc,\Lc}^\times$. Then $G_y=\Sym_f:=\phi^{-1}(\KMon)$.
The partition function of the system $(\Hc(\Dc,\Lc),\sigma_t,\Hc_y,\Ham_y)$,
coincides with the zeta function $\zeta_{\Dc,\Lc}(\beta)$ of $(\Dc,\Lc)$ (see
Definition \ref{zeta-general}).
\end{prop}

Moreover, it follows from \ref{lemme-determinant} that the Hamiltonian
has positive energy in the representation $\pi_y$.

\subsection{Symmetries}
\label{symmetries-general}
Let $(\Dc,\Lc)=((G,X,V,M),(L,K,\KMon))$ be a BCM pair with neat level.
We will denote the center of $G$ by $C$.

Recall that $\Sym_f$ is the semigroup $\phi^{-1}(\KMon)$.
For $m\in \Sym_f$ and $c\in C(\R)$, we define
$$
\theta_{(m,c)}(f)(g,\rho,[z,l]):=
f(g,\rho\phi(m),[cz,l]).
$$

\begin{lemma}
This gives a well defined right action of
$$\Sym(\Dc,\Lc):=\Sym_f(\Dc,\Lc)\times C(\R)$$
on $\Hc(\Dc,\Lc)$ which moreover commutes with the time evolution.
\end{lemma}
\begin{proof}
The action is well-defined because $K$ acts on $Y_{\Dc,\Lc}$ on the
left, while $\Sym$ acts on the right. Recalling that the time evolution
is given by the formula $(\sigma_tf)(g,y)=\bigl(\det \phi(g)\bigr)^{it}f(g,y)$,
it is clear that the action of $\Sym$ commutes with $\sigma_t$.
\end{proof}

Let $C\KMon$ be the center of $\KMon$.

\begin{definition}
Let $\Inn(\Dc,\Lc)$ be the subsemigroup of $\Sym$ defined by
$$
\Inn(\Dc,\Lc):=C(\Q)\cap\phi^{-1}(C\KMon).
$$
\end{definition}

\begin{remark}
There is a (diagonal) inclusions of semigroups
$$\Inn(\Dc,\Lc)\subset \Sym(\Dc,\Lc).$$
This gives a natural action of $\Inn(\Dc,\Lc)$ on
$\Hc(\Dc,\Lc)$.
\end{remark}

\begin{definition}
The semigroup $\Out(\Dc,\Lc):=\Inn(\Dc,\Lc)\backslash\Sym(\Dc,\Lc)$
is called the \emph{outer symmetry semigroup of the BCM system}
$(\Hc(\Dc,\Lc),\sigma_t)$.
\end{definition}

In practical situations, the following hypotheses will often be fulfilled
(see Propositions \ref{compute-out-tot-F} and \ref{compute-out-M2}).
\begin{definition}
\label{bon-centre}
The level structure $\Lc=(L,K,\KMon)$ is called \emph{faithful} if the
image $\phi(C(\Q))$ of the center of $G$ commutes with
$\KMon$, i.e., $\phi(C(\Q))\subset C\KMon$.
The level structure $\Lc$ is called \emph{full} if
the natural morphism
$\Out\to C(\Q)\backslash G(\A_f)$
is surjective; if this morphism is an isomorphism, the $\Lc$ is called
\emph{fully faithful}.
\end{definition}

These symmetries are symmetries up to inner automorphisms.
\begin{prop}
\label{inner-symmetries-general}
There is a morphism
$$\Out(\Dc,\Lc)\to \Out(\Hc(\Dc,\Lc),\sigma_t)$$
to the quotient of the automorphism group of the BCM system
by inner automorphisms of the algebra.
\end{prop}
\begin{proof}
We have to prove that $\Inn$ acts by inner automorphisms.
For $n\in \Inn$, we let $\mu_n$ be
$$
\mu_{n}(g,y)=1\textrm{ if }g\in K.n^{-1},\;
\mu_{n}(g,y)=0\textrm{ if }g\notin K.n^{-1}.
$$
We will show that
$$\theta_{(n,n)}(f)=\mu_nf\mu_n^*,$$
i.e., the action of $\theta_{(n,n)}$ is given by the inner
automorphism corresponding to $\mu_n$.

We have, for all $y\in Y_{\Dc,\Lc}$,
\begin{align*}
(\mu_{n}f\mu_{n}^*)(g,y) & = \sum_{h\in K\backslash G(\A_f),hy\in Y}
\mu_n(gh^{-1},hy)(f\mu_n^*)(h,y),\\
&=\sum_{h\in K\backslash G(\A_f),hy\in Y}
\mu_n(gh^{-1},hy)
\sum_{k\in K\backslash G(\A_f),ky\in Y}
f(hk^{-1},ky)\mu_n^*(k,y),\\
&=\sum_{h,k\in K\backslash G(\A_f),hy,ky\in Y}
\mu_n(gh^{-1},hy)
f(hk^{-1},ky)\mu_n(k^{-1},ky).
\end{align*}
Now, by definition of $\mu_n$, the only nontrivial term of this sum is obtained
when $k^{-1}=n^{-1}$ and $gh^{-1}=n^{-1}$, i.e., $k=n$ and $h=ng$.
Since $n$ is central,
\begin{align*}
(\mu_{n}f\mu_{n}^*)(g,y) & = f(ngn^{-1},ny),\\
 & = f(g,n\rho,[z,ln^{-1}]),\\
 & = f(g,\rho n,[nz,l]),\\
 & = \theta_{(n,n)}(f)(g,y).
\end{align*}
\end{proof}

\section{Comparison with the original Bost-Connes-Marcolli systems}
We want to understand how our systems are related with the usual
Bost-Connes-Marcolli systems in the class number one case.
These class number one systems are called principal BCM systems.
They are directly related to Connes-Marcolli systems defined in
\cite{Connes-Marcolli-I}.

\subsection{Principal BCM systems}
\label{BCMprinc}
Let $(\Dc,\Lc)=((G,X,V,M),(L,K,\KMon))$ be a BCM pair with $(G,X)$ classical.

Let $\Gamma:=G(\Q)\cap K$ and
$$
U^\princ:=\{(g,\rho,z)\in G(\Q)\times \KMon\times X
            \mid g\rho\in \KMon\}.
$$

Let $X^+$ be a connected component of $X$, $G(\Q)^+$ be $G(\Q)\cap G(\R)^+$
(where $G(\R)^+$ is the identity component of $G(\R)$) and $\Gamma_+=G(\Q)^+\cap K$.
Let
$$
U^+:=\{(g,\rho,z)\in G(\Q)^+\times \KMon\times X^+
     \mid g\rho\in \KMon\}.
$$

We have a natural action of $\Gamma^2$ (resp. $\Gamma_+^2$)
on $U^\princ$ (resp. $U^+$) given by
$(g,\rho,z)\mapsto (\gamma_1g\gamma_2^{-1},\gamma_2\rho,\gamma_2z)$.
Let $\Zf^{\princ}_{\Dc,\Lc}$
(resp. $\Zf^+_{\Dc,\Lc}$) be the stacky quotient of $U^\princ$ (resp. $U^+$)
by $\Gamma^2$ (resp. $\Gamma_+^2$).

\begin{definition}
The stack groupoid $\Zf^{\princ}_{\Dc,\Lc}$ is called the
\emph{principal BCM groupoid} for the pair $(\Dc,\Lc)$.
\end{definition}

\begin{prop}
\label{BC+}
Suppose that the natural map $\Gamma\to G(\Q)/G(\Q)^+$ is surjective.
Then the natural map
$$\Zf^+_{\Dc,\Lc}\to \Zf^{\princ}_{\Dc,\Lc}$$
is an isomorphism.
\end{prop}
\begin{proof}
\emph{Surjectivity}: Let $u=(g,\rho,z)\in U^\princ$. We want to show that
there exists $\gamma_1,\gamma_2\in\Gamma$ such that
$(\gamma_1,\gamma_2).u=(\gamma_1g\gamma_2^{-1},\rho,\gamma_2z)\in U^+.$
There exists $\gamma_2\in \Gamma$ with $\gamma_2z\in X^+$ because:
1) the definition (\ref{defshimuradatum}) of a Shimura datum implies that
$\pi_0(X)$ is a $\pi_0(G(\R))$-homogeneous space; and 2)
from our hypothesis and the theorem of real approximation, we get a surjection
$\Gamma\to G(\Q)/G(\Q)^+\cong G(\R)/G(\R)^+.$
Our hypothesis now implies that there exists $\gamma_1\in \Gamma$ such that
$\gamma_1g\gamma_2^{-1}\in G(\Q)^+$. This proves surjectivity.

\emph{Injectivity}:
Now suppose that two points $(g_1,\rho_1,z_1)$ and $(g_2,\rho_2,z_2)$
have the same image in the quotient. Then there exists
$\gamma_1,\gamma_2\in \Gamma$ such that
$(g_1,\rho_1,z_1)=(\gamma_1g_2\gamma_2^{-1},\gamma_2\rho_2,\gamma_2z_2)$.
Since $\gamma_2$ stabilizes $X^+$, it is in $G(\R)^+$, and therefore
also in $\Gamma_+$.
This implies that $\gamma_1$ is in $\Gamma_+$. This proves injectivity.
\end{proof}

We denote by $h(G,K)$ the cardinality of the finite set $G(\Q)\backslash G(\A_f)/K$.

\begin{prop}
\label{principal-full}
If $h(G,K)=1$ then the principal and the full BCM groupoids are the same, i.e.,
the natural map
$$\Zf^{\princ}_{\Dc,\Lc}\to \Zf_{\Dc,\Lc}$$
is an isomorphism.
\end{prop}
\begin{proof}
There is a natural map
$$
\begin{array}{cccc}
\psi: & (\Gamma\backslash G(\Q))\times \KMon\times X & \to &
(K\backslash G(\A_f))\times \KMon\times G(\Q)\backslash (X\times G(\A_f))\\
 & (g,\rho,z) &\mapsto & (g,\rho,[z,1])
\end{array}.
$$
The action of $\gamma_2\in \Gamma$ on the source is given by
$(g,\rho,z)\mapsto (g\gamma_2^{-1},\gamma_2\rho,\gamma_2z)$
and on the range by
$(g,\rho,[z,l])\mapsto (g\gamma_2^{-1},\gamma_2\rho,[z,l\gamma_2^{-1}])$.
Since $\Gamma=K\cap G(\Q)$, we have
$$
\begin{array}{cl}
\psi(\gamma_2\cdot(g,\rho,z)) &=(g\gamma_2^{-1},\gamma_2\rho,[\gamma_2z,1]),\\
                              &=(g\gamma_2^{-1},\gamma_2\rho,[z,\gamma_2^{-1}]),\\
                  &= \gamma_2\cdot\psi(g,\rho,z).
\end{array}
$$
This proves that $\psi$, being $\Gamma$-equivariant, induces a well defined map
$$
\overline{\psi}:
(\Gamma\backslash G(\Q))\underset{\Gamma}{\times} [\KMon\times X] \to
(K\backslash G(\A_f))\underset{K}{\times}
[\KMon\times G(\Q)\backslash (X\times G(\A_f))].
$$

Let us prove that $\overline{\psi}$ is surjective. This will essentially follow
from the equalities $G(\A_f)=K.G(\Q)=G(\Q).K$ (the class number one hypothesis $h(G,K)=1$).

For $(g,\rho,[z,l])\in
(K\backslash G(\A_f))\underset{K}{\times}
[\KMon\times G(\Q)\backslash (X\times G(\A_f))]$, there exists
$\gamma_2\in K$ and $l_2\in G(\Q)$ such that $l=l_2\gamma_2$.
Then, we have the equalities in our quotient space
$$
\begin{array}{ccl}
(g,\rho,[z,l]) & = & (g,\rho,[z,l_2\gamma_2])\\
             & = & (g\gamma_2^{-1},\gamma_2\rho,[z,l_2])\\
         & = &  (g\gamma_2^{-1},\gamma_2\rho,[l_2^{-1}z,1]).\\
\end{array}
$$
There exists $\gamma_1\in K$ and $g_1\in G(\Q)$ such that
$\gamma_1g_1=g\gamma_2^{-1}$ and we have the following equalities in our
quotient space
$$
\begin{array}{ccl}
(g,\rho,[z,l]) & = &  (g\gamma_2^{-1},\gamma_2\rho,[l_2^{-1}z,1])\\
         & = &  (\gamma_1g_1,\gamma_2\rho,[l_2^{-1}z,1])\\
         & = &  \psi(g_1,\gamma_2\rho,l_2^{-1}z).
\end{array}
$$
Thus $\overline{\psi}$ is surjective.

Now we prove that $\overline{\psi}$ is injective.
Suppose that
$$\overline{\psi}(g_1,\rho_1,z_1)=\overline{\psi}(g_2,\rho_2,z_2).$$
Then there exists $\gamma_1\in K$, $\gamma_2\in K$, $\gamma_3\in G(\Q)$
such that
$$
(\gamma_1g_1\gamma_2^{-1},\gamma_2\rho_1,[\gamma_3z_1,\gamma_3\gamma_2^{-1}])
=
(g_2,\rho_2,[z_2,1]).
$$
This implies $\gamma_3=\gamma_2$ and then $\gamma_2\in K\cap G(\Q)=\Gamma$.
But we also have $\gamma_1=g_2\gamma_2g_1^{-1}\in G(\Q)\cap K=\Gamma$.
This shows that
$$
(g_2,\rho_2,z_2)
=
(\gamma_1g_1\gamma_2^{-1},\gamma_2\rho_1,\gamma_2z_1)
$$
with $\gamma_1,\gamma_2\in \Gamma$, i.e., $(g_2,\rho_2,z_2)$ and
$(g_1,\rho_1,z_1)$ are the same
in
$(\Gamma\backslash G(\Q))\underset{\Gamma}{\times} [\KMon\times X]$.
This proves injectivity.

To finish, we prove that the bijection
$\overline{\psi}:\Zf^\princ_{\Dc,Lc}\to \Zf_{\Dc,\Lc}$
is compatible with the groupoid structures.
Let $Y^\princ=\KMon\times X$, and $Y=\KMon\times \Sh(G,X)$.
If $(g,\rho,z)\in \Zf^\princ$, the image of $(\rho,z)\in Y^\princ$ under $g\in G(\Q)$
is given by $(g\rho,gz)\in Y^\princ$. The image of $(\rho,[z,1])\in Y$ under $g$
is given by $(g\rho,[z,g^{-1}])\in Y$, which is equal to $(g\rho,[gz,1])$.
This finishes the proof.
\end{proof}

\begin{definition}
Let $(\Dc,\Lc)$ be a BCM pair with neat level.
The algebra $\Hc_{\princ}(\Dc,\Lc)=C_c(Z^{\princ}_{\Dc,\Lc})$ is called
the \emph{principal BCM algebra} for $(\Dc,\Lc)$.
\end{definition}

\subsection{The Bost-Connes system}
\label{BC-system}
Let $F/\Q$ be a number field. Let $G=\Res_{F/\Q}\G_{m,F}$,
$X_F=G(\R)/G(\R)^+\cong \{\pm 1\}^{\Hom(F,\R)}$, $V=F$ and
$M=\Res_{F/\Q}\Mrm_{1,F}$.
Let $K=\OFh^\times$, $L=\Oc_F$ and $\KMon=\OFh=\Mrm_1(\OFh)$.

\begin{definition}
\label{BC-F}
The pair
$$
\Pc(\Res_{F/\Q}\G_{m,F},X_F)=
(\G_{m,F},X_F,F,\Res_{F/\Q}\Mrm_{1,F}),(\Oc_F,\OFh^\times,\OFh))
$$
is called the
\emph{Bost-Connes pair} for $F$. The corresponding algebra
$\Hc(\Res_{F/\Q}\G_{m,F},X_F)$ is called the \emph{Bost-Connes algebra} for $F$.
\end{definition}

\begin{prop}
In the case $F=\Q$, $\Hc(\G_{m,\Q},\{\pm 1\})$ is the original Bost-Connes algebra.
\end{prop}
\begin{proof}
Recall from Section~\ref{intro-BC-classic} that the Bost-Connes algebra
is the convolution algebra of the groupoid
$Z_{BC}\subset \Q^\times_+\times \Zh$ of pair $(g,\rho)$ with
$g\rho\in\Zh$; thus we need only to show that $Z_{BC}$ coincides with the BCM
groupoid $Z$ of the Bost-Connes pair. Indeed, in the notation of
Section~\ref{BCMprinc}, we have
$$
U^+=\{\,(g,\rho,1)\in \Q^\times_+\times\Zh\times \{1\}\mid g\rho\in\Zh\,\},\quad
\Gamma=\{\pm 1\},\text{ and }\Gamma_+=1.
$$
Therefore $\Z^+:=\Gamma^+\backslash U^+=Z_{BC}$; the map
$\Gamma\to G(\Q)/G(\Q)^+$ is an isomorphism of $\{\pm 1\}$; and
$h(\G_{m,\Q},\Zh^\times)=1$, since it is the usual class number of~$\Q$.
The proposition follows from Propositions~\ref{BC+} and~\ref{principal-full}.
\end{proof}

\subsection{The Connes-Marcolli system}
We now show that in the $\GL_{2,\Q}$ case, we obtain exactly the same
groupoid as Connes and Marcolli \cite{Connes-Marcolli-I}.
This groupoid is only a stack-groupoid, not a usual groupoid.
This restriction was circumvented by Connes and Marcolli using
functions of weight 0 for the scaling action (see \cite{Connes-Marcolli-I},
remark shortly preceeding 1.83). Such a scaling action is not canonically
defined in the general case we consider. As explained before,
we deliberately chose to view this groupoid as a stack-groupoid in order
to define a natural groupoid algebra for it that depends on the resolution
of stack singularities given by the choice of $K$.

Consider the Shimura datum $(\GL_{2,\Q},\Hb^\pm)$, $V=\Q^2$,
and $M=\Mrm_{2,\Q}$.
let $L=\Z^2$, $K=\GL_2(\Zh)$, and $\KMon=\Mrm_2(\Zh)$.

\begin{definition}
The pair
$$
\Pc(\GL_2,\Hb^\pm):=
((\GL_2,\Hb^\pm,\Q^2,\Mrm_{2,\Q}),(\Z^2,\GL_2(\Zh),\Mrm_2(\Zh)))
$$
is called \emph{the modular BCM pair}. The corresponding BCM stack-groupoid is
denoted by $\Zf_{\GL_2,\Hb^\pm}$.
\end{definition}

The stack-groupoid $\Zf^+_{\GL_2,\Hb^\pm}$ is defined as in Section
\ref{BCMprinc}. This is exactly the groupoid studied
by Connes and Marcolli in \cite{Connes-Marcolli-I}.

\begin{lemma}
\label{comparaison-CM}
Our BCM stack-groupoid is the same as Connes and Marcolli's one. In other
words, the natural map
$$
\Zf^+_{\GL_2,\Hb}\to \Zf_{\GL_2,\Hb^\pm}
$$
is an isomorphism.
\end{lemma}
\begin{proof}
We have in this case $h(G,K)=1$ so by Proposition \ref{principal-full},
we have $[Z^{\princ}_{\GL_2,\Hb^\pm}]\cong [Z_{\GL_2,\Hb^\pm}]$.
The map $\GL_2(\Z)\to \GL_2(\R)/\GL_2(\R)^+$ is surjective,
so that we can apply proposition \ref{BC+}, which tells us that
$\Zf^+_{\GL_2,\Hb^{pm}}\cong \Zf^{\princ}_{\GL_2,\Hb^\pm}$.
\end{proof}

\section{Operator theoretic results on BCM algebras}
\subsection{The C*-algebra associated to a BCM datum}
Let $(\Dc,\Lc)=((G,X,V,M),(L,K,\KMon))$ be a BCM pair with neat level.
On the algebra $\Hc(\Dc,\Lc)$, we put the following norm: for every
$f\in\Hc(\Dc,\Lc)$,
$$
\|f\|=\underset{y\in Y_{\Dc,\Lc}}{\sup}\|\pi_y(f)\|.
$$

\begin{lemma}
\label{C*norm}
This defines a C*-norm on $\Hc(\Dc,\Lc)$, i.e., $\|f^*f\|=\|f\|^2$.
\end{lemma}
\begin{proof}
Indeed, it is easy to check that this is a seminorm satisfying the
C*-condition (Definition~\ref{C*-condition}): observe that for
arbitarily small $\epsilon>0$ there is a $y$ such that
$\norm{f}^2-\epsilon=\norm{\pi_y(f)}^2$.  We then have
\[
\norm{f^*f}
\ge\norm{\pi_y(f^*f)}
=\norm{\pi_y(f)}^2=\norm{f}^2-\epsilon,
\]
which of course means that $\norm{f^*f}\ge\norm{f}^2$.  This inequality
is easily shown to imply the C*-condition.

That we get a norm (i.e., $\|f\|=0$ only when $f=0$), and not just a
seminorm, follows from the fact that $f(g,y)\neq 0$ implies that
$\pi_y(f)\neq0$:
\[
\innerprod{\pi_y(f)\varepsilon_g}{\varepsilon_g}
=f(1,gy)
=f(g,y)
\neq0.
\]
Here $\epsilon_g\in \Hc_y$ is the unit vector which takes value $1$ at $g$, and
$0$ elsewhere.
\end{proof}

\begin{definition}
The completion of $\Hc(\Dc,\Lc)$ under the norm $\|.\|$ is denoted
$\Ac(\Dc,\Lc)$ and called the \emph{BCM C*-algebra}.
\end{definition}

\subsection{Construction of extreme KMS states at small temperature}
Let $(\Dc,\Lc)=((G,X,V,M),(L,K,\KMon))$ be a summable BCM pair.
Recall that
$Y_{\Dc,\Lc}^\times=\{(g,\rho,[z,l])\in Y_{\Dc,\Lc}\mid\rho\textrm{ invertible}\}$.

\begin{lemma}
Let $y\in Y_{\Dc,\Lc}^\times$.
Let $\beta$ be such
that the zeta function $\zeta_{\Dc,\Lc}(\beta)$ converges.
The state
$$\Phi_{\beta,y}(f):=\frac{\Trace(\pi_y(f)e^{-\beta\Ham_y})}{\zeta_{\Dc,\Lc}(\beta)}$$
is a $\KMS$
state for the system $(\Ac(\Dc,\Lc),\sigma_t)$ (by Lemma~\ref{lemme-determinant},
$\zeta_{\Dc,\Lc}(\beta)\neq 0$).
\end{lemma}
\begin{proof}
By construction, the algebra $\Hc(\Dc,\Lc)$ is a norm-dense subalgebra of
$\Ac(\Dc,\Lc)$, which is also $\sigma_z$-invariant. Thus, to verify the
$\KMS$ condition, it is enough to show that
$$\Phi_{\beta,y}(f_1\sigma_{i\beta}(f_2))=\Phi_{\beta,y}(f_2f_1)$$
for every pair of functions $f_1,f_2\in \Hc(\Dc,\Lc)$; see Proposition
\ref{KMS-prop}. The convergence of the zeta function implies
that the operator $e^{-\beta H_y}$ is trace class.
The invariance of the trace under cyclic permutations implies that
$$
\begin{array}{ccl}
\zeta_{\Dc,\Lc}(\beta).\Phi_{\beta,y}(f_1\sigma_{i\beta}(f_2)) &=&
\Trace(f_1e^{-\beta H_y}f_2e^{\beta H_y}e^{-\beta H_y}),\\
&=&\Trace(f_1e^{-\beta H_y}f_2),\\
&=&\Trace(f_2f_1e^{-\beta H_y}),\\
&=&\zeta_{\Dc,\Lc}(\beta).\Phi_{\beta,y}(f_2f_1),
\end{array}
$$
which finishes the proof of the KMS condition.
\end{proof}

The commutant of a subset $S\subset \Bc(\Hc_y)$ is by definition
$S'=\{a\in \Bc(\Hc_y)\mid as=sa,\forall s\in S\}$.
\begin{lemma}
\label{irreducible-piy}
If $y\in Y_{\Dc,\Lc}^\times$, then
the commutant $\pi_y(\Ac(\Dc,\Lc))'$ consists only of scalar operators.
\end{lemma}
\begin{proof}
In general, if $y\in Y_{\Dc,\Lc}$, then the Von Neumann algebra
$\pi_y(\Ac)'$ is
generated by the right regular representation of the isotropy group
$Z_{y,y}:=\{[g,y]\in Z\mid s[g,y]=[y]=[gy]=t[g,y]\}$
(cf. \cite{Connes7} Proposition VII.5).
If $y$ is now in $Y_{\Dc,\Lc}^\times$, then the isotropy group
$Z_{y,y}$ is trivial. Therefore, the commutant $\pi_y(\Ac)'$ consists
only of scalar operators.
\end{proof}

Recall that the set of $\KMS$ states is a convex simplex (see Proposition
\ref{KMS-structure}), whose extreme points
are called \emph{extreme $\KMS$ states}.
\begin{prop}
\label{prop-construction-KMS}
Let $y\in Y_{\Dc,\Lc}^\times$ be an invertible element of $Y_{\Dc,\Lc}$. Let $\beta$ be such
that the zeta function $\zeta_{\Dc,\Lc}(\beta)$ converges.
The $\KMS$ state
$$
\Phi_{\beta,y}(f):=
\frac{\Trace(\pi_y(f)e^{-\beta\Ham_y})}{\zeta_{\Dc,\Lc}(\beta)}
$$
is extremal of type $\mathrm{I}_\infty$.
\end{prop}
\begin{proof}
By Proposition \ref{KMS-structure}, the property, for $\Phi_{\beta,y}$,
of being extreme is equivalent to the property of being a factor state, i.e.,
the algebra $\Ac(\Dc,\Lc)$ generates a factor in the GNS representation
of $\Phi_{\beta,y}$.
Following Harari-Leichtnam, \cite{Harari-Leichtnam}, proof of Theorem 5.3.1,
the GNS representation is (up to unitary equivalence)
$$\tilde{\pi}_y=\pi_y\otimes \id_{\Hc_y}:\Ac(\Dc,\Lc)\to \Bc(\Hc_y\otimes\Hc_y),$$
and the associated cyclic vector is
$$
\Omega_{\beta,y}=
\zeta_{\Dc,\Lc}(\beta)^{-1/2}\sum_{h\in K\backslash G_y}
\det(\phi(h))^{-1/2}\epsilon_h\otimes\epsilon_h,
$$
where $\epsilon_h$ is the basis vector of $\Hc_y$ that takes values $1$ at $h$,
and $0$ elsewhere.

The properties that characterize the
triple~$(\Hc_y\otimes\Hc_y,\tilde\pi_y,\Omega_{\beta,y})$ as the
GNS representation of~$\Phi_{\beta,y}$ are precisely:
\begin{enumerate}
\item
$\Phi_{\beta,y}(f)
=\innerprod{\tilde\pi_y(f)\Omega_{\beta,y}}{\Omega_{\beta,y}}$, for every
$f\in \Ac(\Dc,\Lc)$; and
\item The orbit $\tilde\pi_y(\Ac(\Dc,\Lc))\Omega_{\beta,y}$ is dense in the Hilbert
space~$\Hc_y\otimes\Hc_y$.
\end{enumerate}
These two properties are verified by direct calculation.  For example,
to verify the second condition first observe that
\[
\pi_y(f)\epsilon_h
=\sum_{g\in K\backslash G_y}f(gh^{-1},hy)\,\epsilon_g,
\]
and so
\[
\tilde\pi_y(f)\Omega_{\beta,y}
=\zeta_{\Dc,\Lc}(\beta)^{-1/2}\sum_{g,\,h\in K\backslash G_y}
 \det(\phi(h))^{-\beta/2}f(gh^{-1},hy)\,\epsilon_g\otimes\epsilon_h.
\]
But since $G_y=\Sym_f$, every $\det(\phi(h))$ is positive, and
we can choose $f$ to have sufficiently small
support about $(gh^{-1},hy)$ to see that the basis vector
$\epsilon_g\otimes\epsilon_h$ lies in the closure of
$\tilde\pi_y(\Ac(\Dc,\Lc))\xi_{\beta,y}$.

By Lemma \ref{irreducible-piy}, we know that the commutant
$\pi_y(\Ac(\Dc,\Lc))'$ consists of scalar operators.
It is then clear that
$\tilde\pi_y(\Ac(\Dc,\Lc))'
=\pi_y(\Ac(\Dc,\Lc))'\otimes \Bc(\Hc_y)=\C\otimes \Bc(\Hc_y)$, and so
\[
\tilde\pi_y(\Ac(\Dc,\Lc))''
=\Bc(\Hc_y)\otimes\C\id_{\Hc_y}
\cong B(\Hc_y).
\]
This proves that $\Phi_{\beta,y}$ is a Type~$\mathrm{I}_\infty$ factor
state.
\end{proof}

\begin{question}
\label{classif-etats-KMS-general}
Let $(\Dc,\Lc)=((G,X,V,M),(L,K,\KMon))$ be a BCM pair.
Is is true that for $\beta>\!\!>0$,
the map $y\mapsto \Phi_{\beta,y}$ induces a bijection from the Shimura variety
$\Sh(G,X)$
to the space $\Ec_\beta$ of extremal $\KMS$ states on $(\Hc(\Dc,\Lc),\sigma_t)$?
\end{question}

\section{A Bost-Connes system for number fields}
\label{BCM-number-fields}
\subsection{Reminder of Dedekind zeta functions}
\label{Dedekind-zeta}
A first step in understanding what a good analogue of the Bost-Connes algebra
may be is to find a nice description of the partition functions.
This was done first
by Harari and Leichtnam in \cite{Harari-Leichtnam} in the class number one case
and by Paula Cohen in \cite{Cohen} for general number fields, where she used
an adelic description of the Dedekind zeta function.

Let $F$ be a number field.
There is a multiplicative semigroup injection
$\Oc_F^\times\to \Oc_F$.
Let $\OFh^\natural:=\OFh\cap\A_{f,F}^\times$ and
$\Zh^\natural:=\Zh\cap \A_{f,\Q}^\times$.
Then the space $\OFh^\times\backslash\OFh^\natural$ is identified with the multiplicative
semigroup $I_F$ of integral ideals in $F$.
The norm map induces a natural map
$$
\Nm:\OFh^\times\backslash\OFh^\natural\to
\Zh^\times\backslash\Zh^\natural\cong \Z^\times\backslash\Z-\{0\}
\cong \N^\times.
$$
This is the usual norm on ideal classes.

Now the Dedekind zeta function of $F$ can be expressed as
$$\zeta_F(s)=\sum_{n\in \OFh^\times\backslash\OFh^\natural}\frac{1}{\Nm(n)^s}.$$

\subsection{The adelic Bost-Connes algebra}
\label{BC-algebra}
We recall from Subsection \ref{BC-system} the definition of the
Bost-Connes datum for number fields.
Let $F/\Q$ be a number field. Let $G=\Res_{F/\Q}\G_{m,F}$ and
$X_F=G(\R)/G(\R)^+\cong \{\pm 1\}^{\Hom(F,\R)}$.

Following Definition \ref{BC-F}, the Bost-Connes pair for $F$ is
$$
\Pc_F:=\Pc(\Res_{F/\Q}\G_{m,F},X_F)=
((\Res_{F/\Q}\G_{m,F},X_F,F=\End_F(F)),(\Oc_F,\OFh^\times,\OFh)).
$$

Remark that in this case, we have
$$\Sh(\Res_{F/\Q}\G_{m,F},X_F)\cong \pi_0(C_F),$$
where $C_F$ is the idele class group of~$F$.

Then $Y_F=\OFh\times \Sh(\Res_{F/\Q}\G_{m,F},X_F)$ and
$$U_F\subset \A_{f,F}^\times\times Y_F$$
is the subspace of tuples $(g,\rho,[z,l])$ such that $g\rho\in \OFh$.
We let $\gamma_1,\gamma_2\in \OFh^\times\times \OFh^\times$
act on $(g,y=(\rho,[z,l]))\in U_F$ by
$$
(g,\rho,[z,l])
\mapsto
(\gamma_1g\gamma_2^{-1},\gamma_2\rho,[z,l\gamma_2^{-1}]).$$
Let $Z_F:=(\OFh^\times\times \OFh^\times)\backslash U_F$
and let $\Hc_F=\Hc(\Res_{F/\Q}\G_{m,F},X_F):=C_c(Z_F)$ be the corresponding
Bost-Connes algebra for $F$.

\subsection{Partition function and symmetries}
\label{partition-symmetries-BC}
\begin{lemma}
Let $y\in Y_F$ and $\Hc_y=\ell^2(K\backslash G_y)$.
The time evolution on $\Hc_F$ is given by
$$\sigma_t(f)(g,y)=\Nm(g)^{it}f(g,y).$$
The Hamiltonian $H_y$ in $\Hc_y$ is given by
$$(H_y\xi)(g)=\log(\Nm(g))\cdot \xi(g).$$
\end{lemma}
\begin{proof}
Notice that for $a\in F$, the determinant of the $\Q$-linear map $x\mapsto a.x$
on $F$ is the norm $\Nm(a)$.
The lemma then follows from Definitions \ref{tev-general} and \ref{ham-general}.
\end{proof}

\begin{lemma}
The finite symmetry semigroup $\Sym_f(\Res_{F/\Q}\G_{m,F},X_F)$ of $\Pc_F$ is
$\OFh^\natural$. Its zeta function $\zeta_{\Pc_F}$ is the Dedekind
zeta function $\zeta_F$ of $F$.
\end{lemma}
\begin{proof}
The description of the symmetry semigroup follows from its Definition
\ref{sympm-general}. The description of the zeta function follows
from Definition \ref{zeta-general} and Subsection \ref{Dedekind-zeta}.
\end{proof}

We now identify the action of the full symmetry semigroup
$\Sym(\Res_{F/\Q}\G_{m,F},X_F)=\OFh^\natural\times \Res_{F/\Q}\G_{m,F}(\R)$, which
contains archimedean information.

\begin{prop}
\label{compute-out-tot-F}
We have $\Inn(\Res_{F/\Q}\G_{m,F},X_F)=\Oc_F^\natural:=\Oc_F-\{0\}$
and the outer symmetry semigroup
$\Out(\Res_{F/\Q}\G_{m,F},X_F)$ acts on the BCM algebra $\Hc_F$ through
$$\pi_0(F^\times\backslash \A_F^\times).$$
\end{prop}
\begin{proof}
Recall from Eq. \ref{equa-limproj-groupoide} that $\Zf_F$ can be written
as a projective limit
of groupoids $\Zf_{K'}$ for $K'\subset G(\A_f)$ compact open. The $\Sym$-action can
thus be enhanced to an action of the projective limit semigroup
$\limproj_{K'}\Sym/K'$ over all compact open $K'\subset G(\A_f)$.
We know from \cite{De4}, 2.2.3, that
$$
\limproj_{K'} F^\times\backslash
(\A_{f,F}^\times/K'\times \pi_0(\Res_{F/\Q}\G_{m,F}(\R))):=
\Sh(\Res_{F/\Q}\G_{m,F},X_F)
\cong \pi_0(F^\times\backslash\A_F^\times).
$$
It thus remains to prove that the natural map
$$
\Oc_F^\natural\backslash (\OFh^\natural\times \pi_0(\Res_{F/\Q}\G_{m,F}(\R)))
\to
F^\times\backslash(\A_{f,F}^\times\times \pi_0(\Res_{F/\Q}\G_{m,F}(\R)))
$$
is an isomorphism.
The injectivity of this map is clear because
$\Oc_F^\natural:=\Oc_F-\{0\}=\OFh^\natural\cap F^\times$. Since
$F^\times$ acts transitively on $\pi_0(\Res_{F/\Q}\G_{m,F}(\R))$,
to prove surjectivity it suffices to prove surjectivity of the
upper map of the following diagram:
$$
\xymatrix{
\OFh^\natural \ar[r]\ar[d]&
F^\times\backslash\A_{f,F}^\times\ar[d]\\
\Oc_F^\natural\backslash\OFh^\natural/\OFh^\times\ar[r]^\sim&
F^\times\backslash\A_{f,F}^\times/\OFh^\times}
$$
The lower arrow is an isomorphism because these two groups
are equal to the ideal class group of $F$.
Let $g\in \A_{f,F}^\times$ be a finite idele. Then its image by
the vertical projection gives an ideal class, which is the image of
some $m\in \OFh^\natural$. We have $[m]=[g]$ in the right quotient so that
there exists $k\in \OFh^\times$ such that $g=mk \mod F^\times$.
Then $mk\in \OFh^\natural$ is in the preimage of the upper
arrow of the diagram, which proves surjectivity.
\end{proof}

\begin{remark}
Analogous results were already obtained for
$F$ imaginary quadratic by Connes-Marcolli-Ramachandran
(see \cite{Connes-Marcolli-Ramachandran}). In
this case, the datum $(\Res_{F/\Q}\G_{m,F},X_F)$ is classical so that the system
is simpler.
\end{remark}

\section{A Bost-Connes system for Dirichlet characters}
\label{BCM-Dirichlet}
\subsection{Reminder of zeta functions of Dirichlet characters}
\label{Dirichlet-zeta}
We here recall from Neukirch's book \cite{Neukirch}, p. 501,
some facts about characters.

\begin{definition}
A \emph{Hecke character} is a character of the idele class group
$C_F:=\A_F^\times/F^\times$, i.e.,
a continuous homomorphism
$\chi:C_F\to S^1$
to the group $S^1$ of complex numbers of norm $1$.
A \emph{Dirichlet character} is a Hecke character that factors through
the quotient group
$(F_\R)^\times_+\backslash\A_F^\times/F^\times$
where $+$ denotes the connected component for the real topology.
\end{definition}

Let $\mfk=\prod_\pfk \pfk^n$ be a full ideal of $\Oc_F$ and
let $K(\mfk)$ be the kernel of the natural map
$$\OFh^\times\to (\OFh/\mfk)^\times.$$
We say that $\mfk$ is \emph{a module of definition} for the Dirichlet character
$\chi$ if $\chi(K(\mfk))=1$.
We then call $K(\mfk)$ \emph{a subgroup of definition} for $\chi$.

Each Dirichlet character has a module of definition
and for such an $\mfk$, we have a factorisation
$\chi:C(\mfk)\to S^1$,
where $C(\mfk)=((F_\R)^\times_+\times K(\mfk))\backslash\A_F^\times/F^\times$ is
the big ray class group modulo $\mfk$.
Such an $\mfk$ that is moreover minimal (among the modules of definition)
is called \emph{the conductor of the Dirichlet character}.

Recall that $\OFh^\natural=\A_{f,F}^\times\cap \OFh$.
If $\chi:\A_F^\times\to S^1$ is a Dirichlet character, we factor it through
$(F_\R)^\times_+\backslash\A_F^\times$, and thus restrict it to
$\pi_0(F_\R^\times)\times \OFh^\natural$.
Let $K(\mfk)\subset \OFh^\times$ be a primitive
subgroup of definition for $\chi$
and let $K^\natural(\mfk):=\{n\in \OFh^\natural\mid\bar{n}=1\in\OFh/\mfk\}$.

There is an injective map
$K(\mfk)\backslash K^\natural(\mfk)
\hookrightarrow
\OFh^\times\backslash\OFh^\natural$
whose image is the semigroup of all ideals of $F$ prime to $\mfk$.

At least if $\chi$ is trivial at infinity, it induces
$\chi:K(\mfk)\backslash K^\natural(\mfk)\to S^1.$
Now, we can define the $L$-function of our Dirichlet character $\chi$ as
$$
L_F(s,\chi)=
\sum_{n\in K(\mfk)\backslash K^\natural(\mfk)}\frac{\chi(n)}{\Nm(n)^s},
$$
where $\Nm$ was defined in section \ref{Dedekind-zeta}.
In the particular case of a class character, we have
$$
L_F(s,\chi)=
\sum_{n\in \OFh^\times\backslash\OFh^\natural}
\frac{\chi(n)}{\Nm(n)^s}.
$$

\subsection{A Bost-Connes algebra for Dirichlet characters}
Let $\chi:\A_F\to S^1$ be a Dirichlet character that is supposed to be trivial
at infinity.
Let $G=\Res_{F/\Q}\G_{m,F}$ and
$X:=G(\R)/G(\R)^+\cong \{\pm 1\}^{\Hom(F,\R)}$.
Let $\mfk$ be the conductor of $\chi$ and
$\KMon(\mfk)\subset \OFh$ be the multiplicative semigroup defined by
$$
\KMon(\mfk)=
\Ker_{\mathrm{mult}}(\OFh\to\OFh/\mfk):=
\{n\in \OFh\mid\bar{n}=1\in\OFh/\mfk\}.
$$
Recall that we denoted $K(\mfk)\subset \OFh^\times$ the subgroup
$K(\mfk)=\Ker(\OFh^\times\to (\OFh/\mfk)^\times).$
Let $L=\Oc_F$ and $\phi:G\to \GL_\Q(F)$ be the regular representation.

\begin{definition}
The tuple $\Dc_{F,\mfk}:=((\Res_{F/\Q}\G_{m,F},X,K(\mfk)),(\KMon(\mfk),\phi,L))$
is called the \emph{Bost-Connes datum of conductor $\mfk$}.
\end{definition}

The time evolution and Hamiltonian are the same as in the Bost-Connes case
studied in Subsection \ref{partition-symmetries-BC}.

Let $a_\chi$ be the operator on $\Hc_y$ defined by
$$(a_\chi\xi)(g)=\chi(g).\xi(g).$$

\begin{definition}
The \emph{$\chi$-twisted trace} $\Trace_\chi$ on $\Bc(\Hc_y)$ is defined by
$$\Trace_\chi(D)=\Trace(a_\chi.D).$$
\end{definition}

\begin{definition}
The \emph{$\chi$-twisted partition function} of $\Dc_{F,\mfk}$ is defined as
$$\zeta_{\Dc_{F,\mfk},\chi}(s)=\Trace_\chi(e^{-\beta H_y}).$$
\end{definition}

\begin{lemma}
The $\chi$-twisted partition function of $\Dc_{F,\mfk}$ is equal to
the Dirichlet $L$-function $L_F(s,\chi)$.
\end{lemma}
\begin{proof}
This follows from the definition and Subsection \ref{Dirichlet-zeta}.
\end{proof}

Notice that in this case, the symmetry semigroup is not full in the sense
of Definition \ref{bon-centre}.

\begin{remark}
If we want to treat Dirichlet characters with nontrivial infinite component,
it could be useful to construct the groupoid given by the partial action of
$\A_F^\times$ on the space $\A_F\times \pi_0(C_F)$ where
$C_F:=F^\times\backslash\A_F^\times$. If we do the construction as before,
using a quotient by $(\OFh^\times)^2$, the partition function will not be
reasonable. It could be interesting to use Tate's thesis \cite{Tate-these},
that expresses the Dedekind zeta function as an integral, to deal with
this problem. It is not clear to us if a meaningful physical system can be
constructed this way.
\end{remark}

\section{The Hilbert modular BCM system}
We now specialize the general formalism of Section~\ref{BCM-general}
to the case of Hilbert modular Shimura data. This is a good training ground
for the case of a general Shimura datum.

\subsection{Construction}
Let $F$ be a totally real number field.
Let $G:=\Res_{F/\Q}\GL_2$, $X=(\Hb^\pm)^{\Hom(F,\R)}$.
The Shimura datum $(G,X)$ is the \emph{Hilbert modular Shimura datum}.
Let $V$ be the $\Q$-vector space $F^2$ with the natural action $\phi$ of $G$.
Let $M:=\Res_{F/\Q}\Mrm_{2,F}$.
Let $L\subset V$ be $\Oc_F^2$.
Let $K_0=\GL_2(\OFh)\subset G(\A_f)$ and
$\KMon=\Mrm_2(\OFh)\subset \Mrm_2(\A_{f,F})$.
Choose a neat subgroup $K\subset K_0$.

\begin{definition}
The pair $\Pc(G,X,K):=((G,X,V,M),(L,K,\KMon))$ is called
a \emph{Hilbert modular BCM pair} for $F$. The BCM algebra
$\Hc(\Pc)$ is called a \emph{Hilbert modular BCM algebra}.
\end{definition}

\begin{lemma}
\label{principal-full-Hilbert}
If we suppose that $F$ has class number one,
then the natural morphism
$$
\Hc(\GL_{2,F},X,K)\to
\Hc_{\princ}(\GL_{2,F},X,K)$$
from the principal to the full Bost-Connes-Marcolli algebra
is an isomorphism.
\end{lemma}
\begin{proof}
The hypothesis implies (in fact is equivalent to)
$h(G,K)=1$. The result then
follows from proposition \ref{principal-full}.
\end{proof}

We now describe more explicitly the time evolution whose construction
was made in Subsection \ref{tev-hamiltonien-general}.

Let
$C:=\Res_{F/\Q}\G_m$, which is the center of $G=\Res_{F/\Q}\GL_2$.
The natural determinant map $\det:G\to C$ induces
$\det:K\backslash G(\A_f)\to C(\Zh)\backslash C(\A_f)$. The norm
map $\Nm:C\to\G_{m,\Q}$ induces
$$
\Nm:C(\Zh)\backslash C(\A_f)\to
\Zh^\times\backslash \A_f^\times\cong
\Z^\times\backslash \Q^\times\cong \Q_+^\times\subset \R_+^\times.
$$
\begin{lemma}
The time evolution on the Hilbert modular BCM algebra
$\Hc(G,X,K)$ is equal to
$$\sigma_t(f)(g,y)=\Nm(\det(g))^{it}f(g,y).$$
\end{lemma}

\subsection{Symmetries}
We apply the general definitions of Subsection
\ref{symmetries-general} to this case.
We see that
$$\Sym_f=\Mrm_2(\OFh)^\natural:=\GL_2(\A_{f,F})\cap \Mrm_2(\OFh).$$

The center of $G=\Res_{F/\Q}\GL_{2,F}$ is $C=\Res_{F/\Q}\G_m$ and the center
of $\Mrm_2(\OFh)$ is $\OFh$ as a diagonal subsemigroup.
We also have
$\Inn=\Oc_F^\natural:=\Oc_F\cap F^\times$
and an inclusion of semigroups
$\Oc_F^\natural\subset \Mrm_2(\OFh)^\natural.$

The following lemma explains what the symmetries are in the case
of Hilbert modular BCM systems.

\begin{prop}
\label{compute-out-M2}
The outer symmetry semigroup $\Out$ of the Hilbert modular BCM system is
isomorphic to $F^\times\backslash \GL_2(\A_{f,F})\times \Res_{F/\Q}\G_{m,F}(\R)$,
more precisely, the natural map
$$
\Sym_f:=
\Oc_F^\natural\backslash\Mrm_2(\OFh)^\natural\to
F^\times\backslash \GL_2(\A_{f,F})
$$
is an isomorphism.
\end{prop}
\begin{proof}
The injectivity of this map is clear because,
$$
\Oc_F^\natural=F^\times\cap\Mrm_2(\OFh)^\natural.
$$
Let
$(\Mrm_2(\OFh)^\natural)^{-1}=
\{m\in G(\A_f)\mid m^{-1}\in\Mrm_2(\OFh)\}$
be the semigroup of  inverses of elements in $\Mrm_2(\OFh)^\natural$.
We then have
$$\Mrm_2(\OFh)^\natural.(\Mrm_2(\OFh)^\natural)^{-1}=\GL_2(\A_{f,F}).$$
Let $m\in \Oc_F^\natural\backslash\Mrm_2(\OFh)^\natural$. We only
need to prove that $m^{-1}\in\Oc_F^\natural\backslash\Mrm_2(\OFh)^\natural$.
Moreover, to invert a matrix it is enough to prove that its determinant
is invertible.
We have $\det(m)\in \Oc_F^\natural\backslash\OFh^\natural$.
The nonarchimedean part of Proposition \ref{compute-out-tot-F} gives
$\Oc_F^\natural\backslash\OFh^\natural\cong F^\times\backslash \A_{f,F}^\times$,
which implies that
$\det(m)^{-1}\in \Oc_F^\natural\backslash\OFh^\natural
\subset\Oc_F^\natural\backslash\Mrm_2(\OFh)^\natural$.
This finishes the proof.
\end{proof}

\appendix
\section{Stack groupoids}
\label{stack-groupoids}
\subsection{Topological stacks and stack groupoids}
We refer to \cite{Laumon-Moret-Bailly} for the theory of
stacks and to \cite{Noohi} for the theory of topological stacks.

A topological stack will be for us a stack on the site $\Top$ of
topological spaces with usual open coverings, i.e., a category
fibered in groupoids fulfilling some descent condition (which is
precisely described in \cite{Laumon-Moret-Bailly}, D\'efinition
3.1):
\begin{itemize}
\item isomorphisms between two given objects form a sheaf,
\item every descent condition with respect to an open covering is effective.
\end{itemize}
We remark that B. Noohi gave in \cite{Noohi} a more restrictive
condition (saying that the stack admits a covering by a topological
space that is some kind of local fibration), but since we do not use
fine stacky geometry, we will ignore this.

Now we would like to define a stack groupoid. Since both groupoids
and stacks are categories, one needs the language of 2-categories
to describe stack groupoids.
We could do this in a way analogous to the spaces in groupoids
(``espaces en groupo\"ides'') of \cite{Laumon-Moret-Bailly}, 2.4.3.
The theory of Picard stacks, exposed in \cite{SGA4-III}, EXP. XVIII, 1.4
and in \cite{Laumon-Moret-Bailly}, 14.4, is also an inspiring reference.
Our references for the definition
of a weak 2-category are Kapranov-Voevodsky \cite{Kapranov-Voevodsky},
Tamsamani's thesis \cite{Tamsamani}, 1.4,  and Simpson \cite{Simpson}.

Roughtly speaking, a \emph{stack groupoid}
is a groupoid in the category of topological stacks, i.e., the datum
of a tuple $(\Xf_1,\Xf_0,s,t,\epsilon,m)$ composed of two stacks
$\Xf_1$ and $\Xf_0$, equiped with 1-morphisms source $s:\Xf_1\to
\Xf_0$, target $t:\Xf_1\to \Xf_0$, unit $\epsilon:\Xf_0\to
\Xf_1$, and composition
$m:\Xf_1\underset{s,\Xf_0,t}{\times}\Xf_1\to \Xf_1$:
$$
\xymatrix{\Xf_1\ar@<0.5ex>[r]^s \ar@<-0.5ex>[r]_t&
\Xf_0\ar@/_/@<-1ex>[l]_\epsilon},\hspace{2cm}
\Xf_1\underset{s,\Xf_0,t}{\times}\Xf_1\overset{m}{\longrightarrow} \Xf_1.
$$
The 1-morphism
$$
(\Id_{\Xf_1}\times m):
\Xf_1\underset{\Xf_0}{\times} \Xf_1\to
\Xf_1\underset{\Xf_0}{\times} \Xf_1,
$$
that sends morally a pair $(a,b)$ of composable morphisms to the pair
$(a,ab)$,
is supposed to be an equivalence (which implies the existence of an
inverse for the composition law). This tuple should be equiped with
the additional data of an associator
$$\Phi:m\circ (m\times \Id_{\Xf_1}) \overset{\sim}{\Longrightarrow}
         m\circ (\Id_{\Xf_1}\times m),$$
and two unity constraints
$$
U:m\circ (\Id_{\Xf_1}\times \epsilon)\overset{\sim}{\Longrightarrow} \Id_{\Xf_1}\textrm{ and }
V:m\circ (\epsilon\times\Id_{\Xf_1})\overset{\sim}{\Longrightarrow} \Id_{\Xf_1},
$$
fulfilling some higher coherence (or cocycle) conditions:
pentagon, \dots

Rather than writing explicitly the coherence conditions, we prefer
to use Toen's viewpoint of Segal groupoid stacks, which allows one
to forget these conditions by including them in the choice of
inverses for some equivalences in a simplicial diagram.

\begin{definition}
\label{coarse-quotient} Let $(\Xf_1,\Xf_0,s,t,\epsilon,m)$ be a
tuple as before. Its coarse quotient is by definition the quotient
of the coarse moduli space $|\Xf_0|$ (space of isomorphism classes
of objects in $\Xf_0$) by the equivalence relation generated by
$$
x_0\sim x_0'\Leftrightarrow \exists x_1\in |\Xf_1|\textrm{ such that
}s(x_1)=x_0\textrm{ and }t(x_1)=x_0'.$$
\end{definition}

\subsection{Groupoids in the category of spaces with group operations}
Let $\ospace$ be the category of ``spaces with group operation'', i.e.,
pairs $(G,X)$ composed of a
topological space $X$ and a group $G$ that acts on $X$. A morphism
$\phi=(\phi_1,\phi_2):(G_1,X_1)\to (G_2,X_2)$ between two such pairs
is a pair composed of a group morphism $\phi_1:G_1\to G_2$ and a
space morphism $\phi_2:X_1\to X_2$ such that
$$\phi_2(g_1.x_1)=\phi_1(g_1).\phi_2(x_1),\;\forall (g_1,x_1)\in G_1\times X_1.$$

One can define the notion of groupoid in the category $\ospace$.
This is the datum of a tuple $((G_1,X_1),(G_0,X_0),s,t,\epsilon,m)$
fulfilling some natural conditions that we will not write explicitly
here, because we prefer the geometrical language of stacks. There is a
relation between these two languages, which is given by a natural functor
called ``stacky quotient''. Thus one can naturally associate to a groupoid
in $\ospace$ a stack groupoid.

The reason for introducing the category $\ospace$ is to provide an economical
description of the notion of quotient of a groupoid by a group action as
a stack groupoid.

\begin{example}
Let $(\Dc,\Lc)$ be a BCM pair and let
$(U_{\Dc,\Lc},Y_{\Dc,\Lc},s,t,\epsilon,m)$ be the groupoid defined
in Subsection \ref{BCMgroupoid}. There is a natural action of $K^2$
on the groupoid $U_{\Dc,\Lc}$, given by
$$(g,y)\mapsto (\gamma_1g\gamma_2^{-1},\gamma_2y).$$
There is also a natural action of $K$ on $Y_{\Dc,\Lc}$ given by
$$y\mapsto \gamma y.$$
Let $s_K:K^2\to K,(\gamma_1,\gamma_2)\mapsto \gamma_2$ and
$t_K:K^2\to K,(\gamma_1,\gamma_2)\mapsto \gamma_1$ be the two
projections. Then the morphisms in $\ospace$ given by
$(s,s_K),(t,t_K):(K^2,U_{\Dc,\Lc})\to (K,Y_{\Dc,\Lc})$ are called
the equivariant source and target respectively. The fiber product
$$(K^2,U_{\Dc,\Lc})\underset{(s,s_K),(K,Y_{\Dc,\Lc}),(t,t_K)}{\times}(K^2,U_{\Dc,\Lc})$$
is naturally isomorphic to the OSpace
$$(K^2\underset{s_K,K,t_K}{\times}K^2,U_{\Dc,\Lc}\underset{s,Y_{\Dc,\Lc},t}{\times}U_{\Dc,\Lc}).$$
It means that $m$ induces a natural multiplication map
$$m_e=(m_K,m):(K^2,U_{\Dc,\Lc})\underset{(s,s_K),(K,Y_{\Dc,\Lc}),(t,t_K)}{\times}(K^2,U_{\Dc,\Lc})\to
(K^2,U_{\Dc,\Lc})$$ in this equivariant setting. The map
$$m_K:K^2\underset{s_K,K,t_K}{\times}K^2\to K^2$$
is given by
$m(\gamma_1,\gamma_2,\gamma_2,\gamma_3)=(\gamma_1,\gamma_3)$.

Passing to the stacky quotient, we obtain the multiplication map
$$m:\Zf\underset{\Sf}{\times} \Zf\to \Zf$$
on the stack groupoid of Subsection \ref{BCMgroupoid}. Remark that
if the action of $K$ on $Y_{\Dc,\Lc}$ is free, then we obtain a
standard groupoid.
\end{example}

\subsection{Toen's Segal groupoid stacks}
Following the referee's suggestion, we give the general definition
of a stack groupoid. B. Toen proposed us an elegant definition of
the 1-category of stack groupoids which has the advantage of being
very concise. It is based on the simplicial point of view of
2-categories as explained in Tamsamani's thesis \cite{Tamsamani} and
Simpson \cite{Simpson}. Analogous constructions can also be
found in \cite{Toen-Vezzosi-HAG-II}, 1.3.4 and \cite{Leinster}.

Recall that $\Delta$ is the category whose objects are totally
ordered sets $[n]=\{0,\dots,n\}$ and whose morphisms are increasing
maps.

The category of topological stacks can be seen as a 1-category
$\stacks$ with a notion of equivalences.
\begin{definition}
A \emph{Segal stack category} is a simplicial stack
$\Xf_*:\Delta^{op}\to \stacks$ such that the \emph{Segal
morphisms}
$$\Xf_n\to \Xf_1\underset{\Xf_0}{\times}\cdots \underset{\Xf_0}{\times} \Xf_1$$
(given by the $n$ morphisms in $\Delta$, $[1]\to [n]$ that
send $0$ to $i$ and $1$ to $i+1$)
are stack equivalences.
A Segal stack category is a \emph{Segal stack groupoid} if the \emph{right
multiplication morphism}
$$\Xf_2\to \Xf_1\underset{\Xf_0}{\times}\Xf_1$$
given by the two morphisms in $\Delta$
\begin{itemize}
\item $[1]\to[2]$ such that $0\mapsto 0$, $1\mapsto 1$ and
\item $[1]\to[2]$ such that $0\mapsto 0$, $1\mapsto 2$,
\end{itemize}
is a stack equivalence.
\end{definition}

Remark that if we replace the category $\stacks$ by the category
of sets, we find back the usual notions of category and groupoid.

The stacks $\Xf_n$ can be thought as families of $n$ composable
morphisms and the groupoid condition is that the map $(a,b)\mapsto
(a,a\circ b)$ is an equivalence, which implies that each $a$ is an
isomorphism.

The source and target maps $s,t:\Xf_1\to \Xf_0$ are induced by the
morphisms $s=[0]\to [1]:0\mapsto 0$ and $t=[0]\to[1]:0\mapsto 1$.

The choice of an inverse $\phi$ for the Segal morphism
$\Xf_2\to \Xf_1\underset{\Xf_0}{\times}\Xf_1$
allows one to define a \emph{composition}
$\mu:\Xf_1\underset{s,\Xf_0,t}{\times}\Xf_1\to \Xf_2\to \Xf_1$
given by composing $\phi$ with the morphism induced by $[1]\to[2]:0\mapsto 0,
1\mapsto 2$.

The increasing map $[1]\to [0]:0\mapsto 0,1\mapsto 0$ induces a map
$\epsilon:\Xf_0\to \Xf_1$ called the unit map.

Now, choose an inverse $\psi$ to the right multiplication morphism
$$\Xf_2\to \Xf_1\underset{\Xf_0}{\times}\Xf_1,$$
and compose it with
$$
\begin{array}{c}
\Id_{\Xf_1}\times_{\Xf_0}\epsilon:\Xf_1\to \Xf_1\times_{\Xf_0} \Xf_1\quad\textrm{and}\\
d_2:\Xf_2\to \Xf_1,
\end{array}
$$
where $d_2$ is induced by the map $[1]\to [2]:0\mapsto 1,1\mapsto 2$.
Morally, these successive maps send an arrow $a\in \Xf_1$
to the pair $(a,1)$, then to $(a,a^{-1})$ and finally to $a^{-1}$.
Let $i:\Xf_1\to \Xf_1$ be this composition.

To conclude, up to the two additional choices of $\phi$ and $\psi$,
we have obtained a tuple $(\Xf_1,\Xf_0,s,t,\epsilon,i,m)$ giving a
diagram
$$
\xymatrix{\Xf_1\ar@(ul,dl)[]_i\ar@<0.5ex>[r]^s \ar@<-0.5ex>[r]_t&
\Xf_0\ar@/_/@<-1ex>[l]_\epsilon}
$$
and a multiplication $m:\Xf_1\underset{s,\Xf_0,t}{\times}\Xf_1\to
\Xf_1$, which is the basic datum necessary to define any notion of
stack groupoid. The problem is now to define an associativity
2-isomorphism and some other 2-conditions, and to find the right
notion of coherence conditions for them. The point of Toen's
construction is that these coherence
conditions are already encoded in the simplicial
structure. Let's be more explicit.

First remark that since we work in a 2-category, the inverse of a 1-isomorphism
is supposed to be defined up to a \emph{unique} 2-isomorphism.
This implies that the multiplication
$$m:\Xf_1\underset{\Xf_0}{\times}\Xf_1\to \Xf_1$$
is well-defined up to a unique 2-isomorphism.

To define the associator, we use the following strictly commutative
diagrams (products are done over $\Xf_0$)
$$
\xymatrix{
\Xf_3\ar[r]\ar[d] &\Xf_2 \ar[r]\ar[d] & \Xf_1\\
\Xf_2\times \Xf_1\ar[r]\ar[d] &\Xf_1\times\Xf_1\\
(\Xf_1\times\Xf_1)\times\Xf_1 &&
}
$$
and
$$
\xymatrix{
\Xf_3\ar[r]\ar[d] &\Xf_2 \ar[r]\ar[d] & \Xf_1\\
\Xf_1\times \Xf_2\ar[r]\ar[d] &\Xf_1\times\Xf_1\\
\Xf_1\times(\Xf_1\times\Xf_1) &&
}
$$
whose vertical arrows are equivalences.

The uniqueness of inverses of equivalences up to unique
2-isomorphisms gives natural 2-isomorphisms between the
multiplication maps
$$(\Xf_1\times \Xf_1)\times \Xf_1\to \Xf_1$$
and
$$\Xf_1\times (\Xf_1\times \Xf_1)\to \Xf_1$$
and the morphism obtained by choosing an inverse of the equivalence
$$\Xf_3\to \Xf_1\times \Xf_1\times \Xf_1.$$

This gives the associator 2-isomorphism. To check
that the associator fulfils the desired 2-cocycle condition
(pentagon), it is necessary to use the simplicial diagram up to $\Xf_4$.
An explanation of the argument is given in \cite{Leinster}.

\section{Enveloping semigroups}
\label{enveloping-semigroups} In this appendix, we explain what
enveloping semigroups are, as they are a key ingredient in our
formalism (cf. \ref{BCMdata}).

\subsection{Drinfeld's classification}
All groups will be over a field of characteristic $0$.
Recall from Subsection \ref{BCMdata} the following definition.
\begin{definition}
Let $G$ be reductive group over a field. An \emph{enveloping
semigroup} for $G$ is a multiplicative semigroup $M$ such that
$M^\times=G$, $M$ is irreducible and normal.
\end{definition}

Such semigroups were classified by V. Drinfeld in private notes
\cite{Drinfeld-semigroup} that were kindly given to us by L. Lafforgue.
Suppose that the base field is algebraically closed.
Choose a maximal torus $T\subset G$ and a Borel subgroup $B\supset T$. Let $W$
denote the Weyl group for $(G,B,T)$, and let $X=\Hom(T,\G_m)$.
\begin{theorem}[Drinfeld]
There exists a bijection between
\begin{enumerate}
\item the set of normal affine irreducible semigroups $M$ containing $G$ as
their group of units, and
\item the set of $W$-invariant rational polyhedral convex cones
$K\subset X\otimes_\Z\R$ which contain zero and are non-degenerate,
i.e., not contained in a hyperplane.
\end{enumerate}
\end{theorem}

This classification implies that a semisimple group $G$ has
only one enveloping semigroup, namely $G$ itself. This case is for us not
very interesting (because a BCM system with such an enveloping semigroup has a trivial
zeta function) and we would like to construct more interesting semigroups,
in particular, we would like to construct cartesian diagrams
$$
\xymatrix{
G\ar[r]^\phi\ar[d]  & \GL(V)\ar[d]\\
M\ar[r]         & \End(V)
}
$$
for some fixed representation $\phi:G\to\GL(V)$.

For example, for an adjoint Shimura datum $(G,X)$ (i.e., $ZG=\{1\}$),
the triviality of the enveloping
semigroup implies that the BCM systems we construct have a trivial partition function.
It is then interesting to construct another Shimura datum with adjoint
datum $(G,X)$ and such that the enveloping semigroup is no longer trivial.

There is a natural method due to Vinberg to ``enlarge the center''
of a given semisimple simply connected group $G$ in order to have
an enveloping semigroup that is universal in a certain sense. For the
sake of brevity, we do not discuss this construction here.

\subsection{Ramachandran's construction of Chevalley semigroups}
There is another way to construct enveloping semigroups quite explicitly, which was
communicated to us by N. Ramachandran.
The construction of N. Ramachandran uses the following theorem of Chevalley
(see \cite{Springer1}, 5.1).

\begin{theorem}[Chevalley]
Let $G$ be an algebraic group over $\Q$, and let $\phi:G\to \GL(V)$ be a
faithful representation of $G$.
There is a tensor construction $T^{i,j}:=V^{\otimes i}\otimes V^{\vee,\otimes j}$
and a line $D\subset T^{i,j}$ such that $\phi(G)\subset \GL(V)$ is the stabilizer of
this line.
\end{theorem}

\begin{definition}
Let $G$ be an algebraic group over $\Q$,
$\phi:G\to \GL(V)$ be a faithful representation of $G$.
Let $T$ and $D$ be as in Chevalley's theorem. Suppose that $T=V^{\otimes i}$ (resp.
$T=V^{\vee,\otimes i}$) contains no (resp. only) dual tensor factors.
The multiplicative semigroup
$$
\begin{array}{c}
M(G,V,T,D):=\{m\in\End(V)\mid m.D\subset D\}\\
(\textrm{resp. }M(G,V,T,D):=\{m\in\End(V)\mid ^tm.D\subset D\})
\end{array}
$$
is called a \emph{Chevalley enveloping semigroup} for $G$ in $\End(V)$.
\end{definition}

\begin{example}
If $G=\GL_2$ and $V$ is the standard representation, then
$D=\bigwedge^2 V\subset V^{\otimes 2}$ is a line as in Chevalley's theorem,
and $\Mrm_2$ is the corresponding Chevalley enveloping semigroup.
\end{example}

\begin{example}
Let $(V,\psi)$ be a $2g$-dimensional $\Q$-vector space equipped with an
alternating form $\psi\in \bigwedge^2 (V^\vee)$. Then the line
$D=\Q\langle\psi\rangle\subset V^{\vee,\otimes 2}$ is a line as in
Chevalley's theorem for the group $\GSp_{2g}$ with its standard
representation, and the points of the
corresponding semigroup in a commutative $\Q$-algebra $A$ are given by
$$
\MSp_{2g}(A):=\{m\in \End(V)(A)\mid \exists\,\mu(m)\in A,\;
\psi(m.x,m.y)=\mu(m)\psi(x,y),\;\forall\,x,y\in V_A\}.
$$
\end{example}

\bibliographystyle{alpha}
\bibliography{$HOME/travail/fred}

\end{document}